\documentclass[10pt]{amsart}

\usepackage{amsmath, amsthm, amssymb, enumerate}
\usepackage{color}
\usepackage{empheq}

\chardef\coloryes=1 %%%out
\chardef\isitdraft=0 %%%out
\ifnum\isitdraft=1 %%%out
   \input macros0.tex %%%out
   \def\eqref#1{({\ref{#1}})}                %saves writing paranthesis%%%out

\usepackage[color,notcite]{showkeys}%%%out
\definecolor{refkey}{gray}{.3}%%%out

\definecolor{labelkey}{gray}{.3}%%%out

\else%%%out

\fi%%%out

\newcommand{\ip}[2]{\left<#1,#2\right>}

\usepackage{amsaddr}

\begin{document}
%\newcommand{\llabel}{\label}             %\llabel is just a synonim for \label
%\newcommand{\rref}{\ref}                 %\rref is  just a synonim for \ref
%\newcommand{\ccite}{\cite}               %\ccite is just a synonim for
                                         %when ref. to equations
\def\Dg{{D'g}}
\def\ua{u^{\alpha}}
\def\intint{\int\!\!\!\!\int}
\def\intinttext{\int\!\!\!\int}
\def\intintint{\int\!\!\!\!\int\!\!\!\!\int}
\def\intintintint{\int\!\!\!\!\int\!\!\!\!\int\!\!\!\!\int}
\def\ques{{\cor \underline{??????}\cob}}
\def\nto#1{{\coC \footnote{\em \coC #1}}}
\def\fractext#1#2{{#1}/{#2}}
\def\fracsm#1#2{{\textstyle{\frac{#1}{#2}}}}   %smaller version of frac
\def\baru{U}
\def\nnonumber{}
\def\palpha{p_{\alpha}}
\def\valpha{v_{\alpha}}
\def\qalpha{q_{\alpha}}
\def\walpha{w_{\alpha}}
\def\falpha{f_{\alpha}}
\def\dalpha{d_{\alpha}}
\def\galpha{g_{\alpha}}
\def\halpha{h_{\alpha}}
\def\psialpha{\psi_{\alpha}}
\def\psibeta{\psi_{\beta}}
\def\betaalpha{\beta_{\alpha}}
\def\gammaalpha{\gamma_{\alpha}}
\def\Talpha{T}
\def\TTalpha{T_{\alpha}}
\def\TTalphak{T_{\alpha,k}}
\def\falphak{f^{k}_{\alpha}}

\def\cof{\mathop{\rm cf\,}\nolimits}
\def\abar{\overline a}
\def\vbar{\overline v}
\def\qbar{\overline q}
\def\wbar{\overline w}
\def\etabar{\overline\eta}
\def\Om{\Omega}
\def\ep{\epsilon}
\def\al{\alpha}
\def\Ga{\Gamma}
\def\R{\mathbb R}

\def\barB{B}

\newcommand{\bv}{{u}}
\newcommand {\Dn}[1]{\frac{\partial #1  }{\partial \nu}}
\def\mm{m}

%if \isitdraft=0 or \coloryes=0%%%out
\def\cor{{}}%%%out
\def\cog{{}}%%%out
\def\cob{{}}%%%out
\def\coe{{}}%%%out
\def\coA{{}}%%%out
\def\coB{{}}%%%out
\def\coC{{}}%%%out
\def\coD{{}}%%%out
\def\coE{{}}%%%out
\def\coF{{}}%%%out
%\def\MR#1{}%%%out

%%3
%\ifnum\isitdraft=1%%%out
\ifnum\coloryes=1%%%out

%   \definecolor{coloraaaa}{rgb}{0.1,0.2,0.8}%%%out
%   \definecolor{colorbbbb}{rgb}{0.1,0.7,0.1}%%%out
%   \definecolor{colorcccc}{rgb}{0.8,0.3,0.9}%%%out
%   \definecolor{colordddd}{rgb}{0.0,.5,0.0}%%%out
%   \definecolor{coloreeee}{rgb}{0.8,0.3,0.9}%%%out
%   \definecolor{colorffff}{rgb}{0.8,0.3,0.9}%%%out
%   \definecolor{colorgggg}{rgb}{0.5,0.0,0.4}%%%out

% \def\cog{\color{green}}%%%out
 \def\cog{\color{colordddd}}%%%out
 \def\cob{\color{black}}%%%out
 \def\cor{\color{red}}%%%out
 \def\coe{\color{colorgggg}}%%%out

 \def\coA{\color{coloraaaa}}%%%out
 \def\coB{\color{colorbbbb}}%%%out
 \def\coC{\color{colorcccc}}%%%out
 \def\coD{\color{colordddd}}%%%out
 \def\coE{\color{coloreeee}}%%%out
 \def\coF{\color{colorffff}}%%%out
 \def\coG{\color{colorgggg}}%%%out

%%4
\fi%%%out
 \ifnum\isitdraft=1%%%out
   \baselineskip=17pt%%%out
   \input macros.tex%%%out
   \def\blackdot{{\color{red}{\hskip-.0truecm\rule[-1mm]{4mm}{4mm}\hskip.2truecm}}\hskip-.3truecm}%%%out
   \def\bdot{{\coC {\hskip-.0truecm\rule[-1mm]{4mm}{4mm}\hskip.2truecm}}\hskip-.3truecm}%%%out
   \def\purpledot{{\coA{\rule[0mm]{4mm}{4mm}}\cob}}%%%out
   \def\purpledottwo{{\coA{\rule[0mm]{4mm}{4mm}}\cob}}%%%out
   \def\pdot{\purpledot}%%%out
\else%%%out  
   \baselineskip=15pt
   \def\blackdot{{\rule[-3mm]{8mm}{8mm}}}%%%out
   \def\purpledot{{\rule[-3mm]{8mm}{8mm}}}%%%out
   \def\pdot{}%%%out
%%5
\fi%%%out

\def\nts#1{{\hbox{\bf ~#1~}}} %nts=note to self
\def\nts#1{{\cor\hbox{\bf ~#1~}}} %nts=note to self%%%out
\def\ntsf#1{\footnote{\hbox{\bf ~#1~}}} %nts=note to self
\def\ntsf#1{\footnote{\cor\hbox{\bf ~#1~}}} %nts=note to self%%%out
\def\bigline#1{~\\\hskip2truecm~~~~{#1}{#1}{#1}{#1}{#1}{#1}{#1}{#1}{#1}{#1}{#1}{#1}{#1}{#1}{#1}{#1}{#1}{#1}{#1}{#1}{#1}\\}%%%out
\def\biglineb{\bigline{$\downarrow\,$ $\downarrow\,$}}%%%out
\def\biglinem{\bigline{---}}%%%out
\def\biglinee{\bigline{$\uparrow\,$ $\uparrow\,$}}%%%out

\def\mbar{{\overline M}}
\def\tilde{\widetilde}
\newtheorem{Theorem}{Theorem}[section]
\newtheorem{Corollary}[Theorem]{Corollary}
\newtheorem{Definition}[Theorem]{Definition}
\newtheorem{Proposition}[Theorem]{Proposition}
\newtheorem{Lemma}[Theorem]{Lemma}
\newtheorem{Remark}[Theorem]{Remark}
\newtheorem{definition}{Definition}[section]
\def\theequation{\thesection.\arabic{equation}}
\def\endproof{\hfill$\Box$\\}
\def\comma{ {\rm ,\qquad{}} }            %comma in a formula
\def\commaone{ {\rm ,\qquad{}} }         %second comma in a formula
\def\dist{\mathop{\rm dist}\nolimits}    %distance
\def\sgn{\mathop{\rm sgn\,}\nolimits}    %sgn
\def\Tr{\mathop{\rm Tr}\nolimits}    %trace
\def\div{\mathop{\rm div}\nolimits}    %divergence
\def\supp{\mathop{\rm supp}\nolimits}    %divergence
\def\divtwo{\mathop{{\rm div}_2\,}\nolimits}    %two dimensional divergence
\def\re{\mathop{\rm {\mathbb R}e}\nolimits}    %distance
\def\indeq{\qquad{}\!\!\!\!}                     %indentation in formulas
\def\period{.}                           %period in a formula
\def\semicolon{\,;}                      %semicolon in a formula
\newcommand{\cD}{\mathcal{D}}
%**end of header

\title[Velocity-vorticity-{V}oigt model]{Global well-posedness of the velocity-vorticity-{V}oigt model of the 3{D} {N}avier-{S}tokes equations}
\date{February 23, 2018}

\author{Adam~Larios}
% \cortext[cor1]{Corresponding author}
\address{Department of Mathematics, 
              University of Nebraska--Lincoln,
%               203 Avery Hall,
              Lincoln, NE 68588, USA}
\email{alarios@unl.edu}

\author{Yuan~Pei}
\address{Department of Mathematics, 
              University of Nebraska--Lincoln,
%               203 Avery Hall,
              Lincoln, NE 68588, USA}
\email{ypei4@unl.edu}

\author{Leo~Rebholz}
\address{Department of Mathematical Sciences, 
              Clemson University,
              Clemson, SC 29634, USA}
\email{rebholz@clemson.edu}
% \thanks{}

% \author{Adam Larios, Yuan Pei, Leo Rebholz} 
% \maketitle
% \date{}
% \bigskip
% \indent Department of Mathematics\\
% \indent University of Nebraska-Lincoln\\
% \indent Lincoln, NE 68588\\
% \indent e-mail: alarios\char'100unl.edu
% \indent e-mail: ypei4\char'100unl.edu\\
% \indent Department of Mathematical Sciences\\
% \indent Clemson University\\
% \indent Clemson, SC 29634
% \indent email: rebholz@clemson.edu

\begin{abstract}
The velocity-vorticity formulation of the 3D Navier-Stokes equations was recently found to give excellent numerical results for flows with strong rotation.  In this work, we propose a new regularization of the 3D Navier-Stokes equations, which we call the 3D velocity-vorticity-Voigt (VVV) model, with a Voigt regularization term added to momentum equation in velocity-vorticity form, but with no regularizing term in the vorticity equation. We prove global well-posedness and regularity of this model under periodic boundary conditions.  We prove convergence of the model's velocity and vorticity to their counterparts in the 3D Navier-Stokes equations as the Voigt modeling parameter tends to zero.   We prove that the curl of the model's velocity converges to the model vorticity (which is solved for directly), as the Voigt modeling parameter tends to zero.  Finally, we provide a criterion for finite-time blow-up of the 3D Navier-Stokes equations based on this inviscid regularization.
\end{abstract}

\maketitle

\noindent\thanks{\em Keywords:\/}
Vorticity-Velocity formulation,
Euler-Voigt,
Navier-Stokes-Voigt,
Global existence,
Inviscid-regularization,
Turbulence models,
Blow-up criteria,
Voigt-regularization,
Turbulence models,
$\alpha$-models,

\noindent\thanks{\em Mathematics Subject Classification\/}:
%
% 35-04, % Explicit machine computation and programs (not the theory of computation or programming)
35A01, % Existence problems: global existence, local existence, non-existence
% 35A02, % Uniqueness problems: global uniqueness, local uniqueness, non-uniqueness
% 35B30, % Dependence of solutions on initial and boundary data, parameters
% 35B41, % Attractors (PDEs)
35B44, % Blow-up
35B65, % Smoothness and regularity of solutions
% 35K51, % Initial--boundary value problems for second-order parabolic systems
% 35K55, % Nonlinear parabolic equations
35Q30, % Navier--Stokes equations
% 35Q31, % Euler equations
35Q35, % PDEs in connection with fluid mechanics
% 35Q86, % PDEs in connection with geophysics 
% 37L30, % Attractors and their dimensions, Lyapunov exponents 
% 65N99, % General Numerical Analysis of BVPs for PDEs
% 76A10, % Viscoelastic fluids
% 76B03, % Existence, uniqueness, and regularity theory (Incompressible inviscid fluids)
76D03, % Existence, uniqueness, and regularity theory (Incompressible viscous fluids)
76D05, % Navier-Stokes equations
% 76D09, % Viscous-inviscid interaction
76D17, %  Viscous vortex flows ( Incompressible viscous fluids )
% 76F05, % Isotropic turbulence; homogeneous turbulence 
% 76F20, % Dynamical systems approach to turbulence [See also 37-XX] 
% 76F25, % Turbulent transport, mixing
% 76F55, % Statistical turbulence modeling [See also 76M35]
% 76F65, % Direct numerical and large eddy simulation of turbulence
76N10 % Existence, uniqueness, and regularity theory
% 76W05, % Magnetohydrodynamics and electrohydrodynamics

\section{Introduction}
\label{sec1}
In recent years, the Voigt-regularization and the velocity-vorticity formulation have seen much study as promising approaches to alleviating some of the analytical and computational difficulty inherent in the 3D Navier-Stokes equations of incompressible fluid flow.  However, as one might expect, neither of these approaches overcomes every difficulty in the equations.  For instance, the Voigt-regularization has a strong regularizing effect, so much so that it destroys certain fundamental qualities of the equations, such as parabolicity and viscosity-driven energy decay.  On the other hand, the velocity-vorticity formulation is merely a reformulation of the equations, and therefore it has no regularizing effect at all, although it is the basis of many well-behaved numerical algorithms.  In this paper, we combine these two approaches, with the intent that the resulting system will retain the best qualities of both systems.  Namely, the intent is that the new system will have solutions that are closer to the actual physics of fluids, while still having enough regularization that the equations are better behaved from the standpoints of mathematical analysis, numerical stability, and computational efficiency.  In this work, we only address the global well-posedness and convergence properties of the system, but a follow-up work will study the numerical and computational properties of the system. 

The incompressible, constant density, 3D Navier-Stokes equations are given by
\begin{subequations}\label{NSE}
\begin{empheq}[left=\empheqlbrace]{align}    
\label{NSE_u}
           &\frac{\partial \tilde{u}}{\partial t} 
           -
           \nu\Delta\tilde{u}
           + 
           (\tilde{u}\cdot\nabla)\tilde{u} 
           + 
           \nabla \tilde{p} 
           = 
           f,
           \\&
           \nabla \cdot \tilde{u} = 0,
           \\&
           \tilde{u}(\cdot, 0) = \tilde{u}_0, 
\end{empheq}
\end{subequations}
 \cite{Constantin_Foias_1988, Temam_2001_Th_Num} for more details on 3D Navier-Stokes equations. 
Note that for smooth solutions, the momentum equation \eqref{NSE_u} can also be written as 
\begin{equation}
\label{NSE_vor}
          \frac{\partial \tilde{u}}{\partial t} 
           -
           \nu\Delta\tilde{u}
           + 
           (\nabla\times\tilde{u})\times\tilde{u} 
           + 
           \nabla(\tilde{p}+\frac{1}{2}|\tilde{u}^2|) 
           = 
           f,
\end{equation}
where $\tilde{u}$ represents the velocity of the fluid, $\tilde{p}$ represents the (density normalized) pressure, and $f$ represents a body force.  
We now propose the following system, which we refer to as the velocity-vorticity-Voigt (or ``VVV'') equations over the three-dimensional periodic box $\mathbb{T}^3 = \mathbb{R}^3/\mathbb{Z}^3=[0, 1]^3$, 
\begin{subequations}\label{Sys1}
\begin{empheq}[left=\empheqlbrace]{align}    
\label{Sys1_u}
           &(I - \alpha^2\Delta)\frac{\partial u}{\partial t} 
           -
           \nu\Delta u
           + 
           w \times u 
           + 
           \nabla p 
           = 
           f,
           \\\label{Sys1_w}&
           \frac{\partial w}{\partial t} 
           - 
           \nu\Delta w 
           + 
           (u\cdot\nabla) w
           - 
           (w\cdot\nabla) u
           =
           \nabla\times f,
           \\&
           \nabla \cdot u = 0,
           \\&
           u(\cdot, 0) = u_0
           \quad\text{ and }\quad
           w(\cdot, 0) = w_0.
\end{empheq}
\end{subequations}
where $u = (u_1, u_2, u_3)$ represents an averaged velocity, $w = (w_1, w_2, w_3)$, which plays the role of vorticity but for which we do not assume $w=\nabla\times u$, and $f$ is an external forcing term. Without loss of generality, we assume for our analysis in later sections that the viscosity $\nu=1$.  Note that in the case where $\alpha=0$, the system formally reduces the velocity-vorticty formulation, while for $\alpha>0$, if one imposes $w = \nabla\times u$, the system formally reduces to the Navier-Stokes-Voigt equations.  

The term $-\alpha^2\Delta\partial_t u$ in \eqref{Sys1_u} is often referred to as the ``Voigt-term'', due to an application of modeling Kelvin-Voigt fluids by A.P. Oskolkov \cite{Oskolkov_1973, Oskolkov_1982} (see also \cite{Kalantarov_1986}).  In the context of the velocity formulation  use of the Voigt term was first proposed as a regularization for either the Navier-Stokes (for $\nu>0$) or Euler (for $\nu=0$) equations in \cite{Cao_Lunasin_Titi_2006}, for small values of the regularization parameter $\alpha$.  This paper also proved global well-posedness of the Voigt-regularized versions of the 3D Euler and 3D Navier-Stokes equations.  These equations have been studied analytically and extended in a wide variety of contexts (see, e.g., \cite{Bohm_1992,Catania_2009,Catania_Secchi_2009,Cao_Lunasin_Titi_2006,Larios_Titi_2009,Larios_Lunasin_Titi_2015,Ebrahimi_Holst_Lunasin_2012,Levant_Ramos_Titi_2009,Khouider_Titi_2008,Olson_Titi_2007,Oskolkov_1973,Oskolkov_1982,Ramos_Titi_2010,Kalantarov_Levant_Titi_2009,Kalantarov_Titi_2009,Kuberry_Larios_Rebholz_Wilson_2012,DiMolfetta_Krstlulovic_Brachet_2015,Layton_Rebholz_2013_Voigt,Larios_Petersen_Titi_Wingate_2015}, and the references therein).  Voigt-regularizations of parabolic equations are a special case of pseudoparabolic equations, that is, equations of the form $Mu_t+Nu=f$, where $M$ and $N$ are (possibly non-linear, or even non-local) operators.  For more about pseudoparabolic equations, see, e.g., \cite{DiBenedetto_Showalter_1981,Peszynska_Showalter_Yi_2009,Showalter_1975_nonlin,Showalter_1975_Sobolev2,Showalter_1972_rep,Carroll_Showalter_1976,Showalter_1970_SG,Showalter_1970_odd,Bohm_1992}.

Directly computing for the vorticity variable has recently become popular in simulations of the incompressible Navier-Stokes equations, as it can be the primary variable of interest in vortex dominated and rotating flows \cite{WB02,LYM06,MF00,G91,GHH90,WWW95,HOR17}.  Such formulations use the equation of vorticity dynamics,
\begin{equation}\label{vort}
\frac{\partial {\tilde w}}{\partial t} - \nu \Delta {\tilde w} +
(\tilde u \cdot \nabla) {\tilde w}  - (\tilde w \cdot\nabla) \tilde u =
 {\nabla \times}\,f,
\end{equation}
and close the system with some relation of $\tilde u$ and $\tilde w$, the most common being $-\Delta \tilde u=\nabla \times \tilde w$, and boundary conditions such as $\tilde w\left.\right|_{\partial\Omega}=\nabla \times \tilde u$.  While such a boundary condition is easily and accurately implementable in finite difference methods on uniform grids, it is not generally appropriate for finite element methods on unstructured meshes.  For this reason, in the recent works \cite{OR10,LOR11,GHOR15,HOR17}, the system is instead closed by coupling to a momentum equation using the $\tilde w$ variable, such as
\[
\tilde u_t + \tilde w\times \tilde u + \nabla \tilde P - \Delta \tilde u = f,
\]
where $\tilde P$ represents the Bernoulli pressure. 
This formulation is able to produce efficient and accurate numerical methods in such settings, in particular due to its use of a natural vorticity boundary condition corresponding to no-slip velocity derived in \cite{GHOR15}.

Because of the numerical successes of such velocity-vorticity systems, it is both natural and important to consider their analysis at the PDE level, as fundamental questions such as well-posedness should be addressed. It is easy to show that determining the global well-posedness of such a system (i.e. \eqref{Sys1} with $\alpha=0$) would solve the Millennium Prize Problem for the 3D Navier-Stokes equations. 
We choose in this work to consider the velocity-vorticity system with a Voigt modeling term in \eqref{Sys1}, for several reasons: first, it allows for analysis of the system to be performed; second, the VVV system limiting behavior as $\alpha\rightarrow 0$ can give insight into the behavior of the $\alpha=0$ case; third, the discretized Voigt term corresponds to a commonly used numerical stabilization for second (and lower) order methods \cite{P06,WL97,DP09,ALP04,LLMNR09}, and thus the VVV system is in this sense the PDE generalization of stabilized numerical discretizations of Navier-Stokes equations in velocity-vorticity form. We note that a steady velocity-vorticity system without regularization terms was analyzed in \cite{ORS17} in the case of no slip velocity boundary conditions, no penetration vorticity boundary conditions, and a natural tangential condition for vorticity that is weakly implemented in a boundary functional involving the pressure; well-posedness of this system was proven, however, it (seemingly) required some new analytic techniques.

We note that we do not add a Voigt modeling term $-\alpha^2\Delta  \frac{\partial w}{\partial t}$ to the vorticity equation of the system.  This could be done, and could make sense as a continuous level generalization of a vorticity equation stabilization.  However, from an analysis point of view, it is more challenging to consider the system \eqref{Sys1} where the Voigt modeling is only applied to the momentum equation, and extension is straightforward for the case when Voigt modeling is also applied to the VVV vorticity equation. Moreover, it may lead to more accurate capturing of the vorticity at the numerical level.  The reason for this is that in rooted in a computational study \cite{Kuberry_Larios_Rebholz_Wilson_2012} of the magnetohydrodynamic (MHD) equations with Voigt-regularization.  Voigt-regularization for MHD was first proposed and studied in \cite{Larios_Titi_2009}, with further study in \cite{Catania_2009,Catania_Secchi_2009,Larios_Titi_2010_MHD}.  The MHD system, with Voigt-regularization added only to the momentum equation, is strikingly similar to system \eqref{Sys1}.    Indeed, if one were to add the term $(w\cdot\nabla)w$ to the right-hand side of equation \eqref{Sys1_u}, the systems would be identical in the $f=0$ case.  In \cite{Kuberry_Larios_Rebholz_Wilson_2012}, it was found that in computational tests of the 2D case on a coarse mesh, putting a Voigt term only on (the MHD analogue of) equation \eqref{Sys1_w} resulted in a better match of level curves of the current density (the analogue of $\nabla\times w$) in fine-mesh simulations than putting a Voigt regularization on both equation, or neither equation.  This is the reason for us only applying Voigt-regularization to equation \eqref{Sys1_u}.

\begin{Remark}
We note that the analysis of the 3D VVV system is somewhat distinct that the analysis of the 3D MHD system with Voigt-regularization added only to the momentum equation.  This is because the cancellation of the nonlinear terms and that occurs in energy estimates for the MHD-Voigt equations does not occur in the VVV system, and therefore one must deal directly with the analogue of the vortex-stretching term $(w\cdot\nabla) u$.  The key is to notice that one may first obtain an energy estimate purely in terms of $u$, and then use this bound to obtain a bound on $w$.  Higher-order estimates on $u$ are not $w$-independent, but can be obtained using a bootstrapping technique, going back and forth between the two equations.
\end{Remark}

\begin{Remark}
 \label{remark_L2V_bdd}
 One might ask whether, in the inviscid (i.e., $\nu=0$) case, results analogous to those in this paper still hold.  This is especially in light of global well-posedness results for the so-called Euler-Voigt equations \cite{Cao_Lunasin_Titi_2006,Larios_Titi_2009}, which are formally the inviscid version of the Navier-Stokes-Voigt equations.  Two fundamental differences arise.  Firstly, the vorticity stretching term $(w\cdot\nabla) u$ can no longer be controlled in the same way, as higher-order derivatives cannot be absorbed into the viscosity.  Thus, one must resort to higher-order estimates, but as in case of the 3D Navier-Stokes equations and related $\alpha$-models \cite{Foias_Holm_Titi_2002,Ilyin_Lunasin_Titi_2006, Chen_Foias_Holm_Olson_Titi_Wynne_1998_PF,Chen_Foias_Holm_Olson_Titi_Wynne_1999,Holm_Titi_2005, Chen_Foias_Holm_Olson_Titi_Wynne_1998_PRL,Cheskidov_Holm_Olson_Titi_2005}, it is far from clear how to close these estimates.  Secondly, in the proof of convergence as $\alpha\rightarrow 0$ (Theorem \ref{T4} below), the estimates depend crucially on the fact that $\int_0^T\|\nabla u(t)\|_{L^2}^2\,dt$ is bounded independently of $\alpha\in(0,1]$, which is a property that one does not have in the Euler-Voigt equations.  Thus, it is not clear how to extend the results of this paper to the inviscid version of the VVV system.
\end{Remark}

The paper is organized as follows. We first provide the necessary preliminaries for our work in subsequent sections in Section~\ref{sec2}, then, we define weak and strong solutions to system \eqref{Sys1}, and state our main theorems. In Section~\ref{sec3}, we prove the existence and uniqueness of global weak solution for \eqref{Sys1} by Galerkin approximation following the ideas from \cite{Constantin_Foias_1988, Temam_2001_Th_Num} (also c.f.\cite{Larios_Pei_triple} with similar approach in full details). In view of the similar behavior of $w$ and the vorticity $\omega = \nabla\times u$, we prove that $w$ indeed tends to $\omega$ in $L^2$ norm as $\alpha\rightarrow 0$  in Section~\ref{sec4}, as well as the convergence of the velocity in \eqref{Sys1} to that of the Navier-Stokes equations. We point out that the numerical and computational studies of the VVV system will be the subject of a forthcoming work.

\section{Preliminaries and Main Results}
\label{sec2}
\subsection{Preliminaries}
\label{subsec2-1}
All through this paper $C$ represents some absolute constant varying line by line, 
and similarly $C_{\alpha}$ indicates the dependence of the constant on $\alpha$. 
We denote $\phi_{j} = \partial \phi/ \partial x_{j}$ 
and 
$\phi_{t} = \partial \phi/\partial t$.  
Also, we denote the mean-free versions of the usual Lebesgue and Sobolev spaces on $\mathbb{T}$ by $L^{p}$ for $1\leq p\leq \infty$ and $H^{s}\equiv W^{s, 2}$ for $s > 0$, respectively; we denote by $C_{w}(I;X)$ the space of weakly continuous functions from an interval $I$ to a Banach space $X$. Let $\mathcal{F}$ be the set of all trigonometric polynomials over $\mathbb{T}^3$ 
and define the subset of $\mathcal{F}$ with divergence-free and zero-average trigonometric polynomials 
$$\mathcal{V} := \left\{ \phi\in\mathcal{F}: \nabla\cdot\phi = 0, \text{ and }\int_{\mathbb{T}^3}\phi\,dx= 0\right\}.$$
We follow the standard convention of denoting by $H$ and $V$ the closures of $\mathcal{V}$ in $L^2$ and $H^1$, respectively, 
with inner products 
$$(v, \tilde{v}) = \sum_{i=1}^3\int_{\mathbb{T}^3}v_{i}\tilde{v}_{i}\,dx \text{ \,\,and\,\,  } ((v, \tilde{v})) = \sum_{i, j=1}^3\int_{\mathbb{T}^3}\partial_{j}v_{i}\partial_{j}\tilde{v}_{i}\,dx,$$
respectively, associated with the norms $\Vert v\Vert_{H}=(v, v)^{1/2}$ and $\Vert v \Vert_{V}=((v, v))^{1/2}$. For the sake of convenience, we use $\Vert v\Vert_{L^2}$ and $\Vert v\Vert_{H^1}$ to denote the above norms in $H$ and $V$, respectively. 
The latter is a norm due to the Poincar\'e inequality 
\begin{equation}
   \label{Poincare}
     \sqrt{\lambda_1}\Vert\phi\Vert_{L^2} \leq \Vert\nabla\phi\Vert_{L^2},
\end{equation}
holding for all  $\phi\in V$, where $\lambda_1$ is the first eigenvalue of the Stokes operator $A$ discussed below. 
Note that we also have the following compact embeddings (see, e.g., \cite{Constantin_Foias_1988, Temam_2001_Th_Num})
$$V \hookrightarrow H_{\sigma} \hookrightarrow V',$$
where $V'$ denotes the dual space of $V$. 

We denote by $H_{\text{curl}}$ a subspace of $H$ whose elements are in $H$ and their curl (taken in the distributional sense) is in $L^2$, i.e., 
$$H_{\text{curl}} := \{ f\in H \left.\right| \nabla\times f \in L^2\},\quad\text{with norm}\quad \Vert f\Vert_{H_{\text{curl}}}:=(\Vert f\Vert_{L^2}^2 + \Vert \nabla\times f\Vert_{L^2}^2)^{1/2},$$ 
and $H_{\text{curl}}^s$ the subspace of $V$ whose elements are in $H^{s}\cap V$ and their curl is in $H^{s}$, i.e., 
$$H_{\text{curl}}^{s} := \{ f\in V \left.\right| \nabla\times f \in H^s\cap V\}\quad\text{with norm}\quad \Vert f\Vert_{H_{\text{curl}}^{s}}:=(\Vert f\Vert_{H^{s}}^2 + \Vert \nabla\times f\Vert_{H^{s}}^2)^{1/2}.$$
For more discussion on the curl-spaces, we refer the readers to \cite{ORS17,GR86} and the references therein. 

The following interpolation result is frequently used in this paper (see, e.g., \cite{Nirenberg} for a detailed proof).
Assume $1 \leq q, r \leq \infty$, and $0<\gamma<1$.  
For $v\in L^q(\mathbb{T}^{n})$, such that  $\partial^\alpha v\in L^{r} (\mathbb{T}^{n})$, for $|\alpha|=m$, then 
\begin{align}\label{PT1}
\Vert\partial_{s}v\Vert_{L^{p}} \leq C\Vert\partial^{\alpha}v\Vert_{L^{r}}^{\gamma}\Vert v\Vert_{L^{q}}^{1-\gamma},
\quad\text{where}\quad
\frac{1}{p} - \frac{s}{n} = \left(\frac{1}{r} - \frac{m}{n}\right) \gamma+ \frac{1}{q}(1-\gamma).
\end{align}

The following results are standard in the study of fluid dynamics, 
in particular for the Navier-Stokes equations and related PDEs,  
and we refer to reader to \cite{Constantin_Foias_1988, Temam_2001_Th_Num} for more details. 
We define the Stokes operator $A:= -P_{\sigma}\Delta$ 
with domain $\mathcal{D}(A):= D(A)$, 
where $P_{\sigma}$ is the Leray-Helmholtz projection. 
Notice that due to the periodic boundary conditions, 
it holds that $A = -\Delta P_{\sigma}$.
Moreover, the Stokes operator can be extended 
as a linear operator from $V$ to $V'$ such that
$$\left<Av, \tilde{v}\right> = ((v, \tilde{v})) \text{  for all  } v,\,\,\tilde{v}\in V.$$
It is well-known that $A^{-1} : H \hookrightarrow \mathcal{D}(A)$ 
is a positive-definite, self-adjoint, compact operator from $H$ into itself,
and $H$ possesses an orthonormal basis of eigenfunctions $\{ w_{k}\}_{k=1}^{\infty}$ of $A^{-1}$, corresponding to a sequence of non-increasing sequence of positive eigenvalues. 
Therefore, $A$ has non-decreasing eigenvalues $\lambda_{k}$, 
i.e., $0 < \lambda_1 \leq \lambda_2, \ldots$ 
since $\{w_{k}\}_{k=1}^{\infty}$ are also eigenfunctions of $A$. 
Furthermore, for any integer $M > 0$, we define $H_{M}:= \text{span}\{w_1, w_2, \ldots, w_{M}\}$ 
and $P_{M} : H \to H_{M}$ be the $L^2$-orthogonal projection onto $H_{M}$. 
Next, for any $v, \tilde{v}, w \in \mathcal{V}$, 
we introduce the convenient notation for the bilinear term
\begin{align*}
     B(v, \tilde{v}) := P_{\sigma}((v\cdot\nabla)\tilde{v}),
\end{align*}
which can be extended to a continuous map 
$B : V \times V \to V'$, such that, for smooth functions $v,\tilde{v},w\in\mathcal{V}$,
\begin{align*}
     \left<B(v, \tilde{v}), w\right> = \int_{\mathbb{T}^3}(v\cdot\nabla) \tilde{v}\cdot w\,dx.
\end{align*} 
Moreover ${\barB}$ has certain symmetry properties, and can be extended as a continuous map ${\barB}$ on various spaces, as in the following lemma, which is proved in, e.g., \cite{Constantin_Foias_1988, Foias_Manley_Rosa_Temam_2001}.
\begin{Lemma}
\label{L1}
The operator $B$ can be extended to an operator, still denoted by $B$, over spaces indicated below, with the following properties. 
\begin{subequations}
\begin{align}
    \label{embd1}
    \ip{B(u,v)}{w}_{V'} &= -\ip{B(u,w)}{v}_{V'}, 
    &\quad\forall\;u\in V, v\in V, w\in V,\\
    \label{embd2}
    \ip{B(u,v)}{v}_{V'} &= 0,
    &\quad\forall\;u\in V, v\in V, w\in V,\\
    \label{B:326}
    |\ip{B(u,v)}{w}_{V'}|
    &\leq C\Vert u\Vert_{L^2}^{1/2} \Vert\nabla u\Vert_{L^2}^{1/2} \Vert\nabla v\Vert_{L^2} \Vert\nabla w\Vert_{L^2},
    &\quad\forall\;u\in V, v\in V, w\in V,\\
    \label{B:623}
    |\ip{B(u,v)}{w}_{V'}|
    &\leq C\Vert\nabla u\Vert_{L^2} \Vert\nabla v\Vert_{L^2} \Vert w\Vert_{L^2}^{1/2} \Vert\nabla w\Vert_{L^2}^{1/2},
    &\quad\forall\;u\in V, v\in V, w\in V,\\
    \label{B:236}
    |\ip{B(u,v)}{w}_{V'}|
    &\leq C\Vert u\Vert_{L^2} \Vert\nabla v\Vert_{L^2}^{1/2} \Vert Av\Vert_{L^2}^{1/2} \Vert\nabla w\Vert_{L^2},
    &\quad\forall\;u\in H, v\in \mathcal{D}(A), w\in V,\\
    \label{B:632}
    |\ip{B(u,v)}{w}_{V'}|
    &\leq C\Vert\nabla u\Vert_{L^2} \Vert\nabla v\Vert_{L^2}^{1/2} \Vert Av\Vert_{L^2}^{1/2} \Vert w\Vert_{L^2},
    &\quad\forall\;u\in V, v\in \mathcal{D}(A), w\in H,\\
    \label{B:i22}
    |\ip{B(u,v)}{w}_{V'}|
    &\leq C\Vert\nabla u\Vert_{L^2}^{1/2} \Vert Au\Vert_{L^2}^{1/2} \Vert\nabla v\Vert_{L^2} \Vert w\Vert_{L^2},
    &\quad\forall\;u\in \mathcal{D}(A), v\in V, w\in H,\\
    \label{B:623s}
   |\ip{B(u,v)}{w}_{V'}|
   &\leq C\Vert\nabla u\Vert_{L^2} \Vert \nabla w\Vert_{L^2} 
   \Vert v\Vert_{L^2}^{1/2}\Vert \nabla v\Vert_{L^2}^{1/2},
   &\quad \forall\;u\in V, v\in V, w\in V,\\
   \label{B:263}
    |\ip{B(u,v)}{w}_{V'}|
    &\leq C\Vert u\Vert_{L^2} \Vert Av\Vert_{L^2} \Vert w\Vert_{L^2}^{1/2} \Vert\nabla w\Vert_{L^2}^{1/2},
    &\quad\forall\;u\in H, v\in \mathcal{D}(A), w\in V.
\end{align}
\end{subequations}
\end{Lemma}

If we formally apply $P_\sigma$ to equation \eqref{Sys1_u} and \eqref{Sys1_w}, we obtain the following functional formulation of system \eqref{Sys1}
\begin{subequations}\label{Sys1_fun}
 \begin{empheq}[left=\empheqlbrace]{align}
   \label{Sys1_fun_u}
     &\frac{d}{dt}(u + \alpha^2 Au)
     +
     Au
     +
     P_\sigma(w\times u)
     =
     P_\sigma f, 
     \\
     \label{Sys1_fun_w}
     &\frac{d w}{dt}
     +
     A w
     +
     {\barB}(u, w)
     -
     {\barB}(w, u)
     =
     \nabla\times f,
   \end{empheq}
\end{subequations}
Formulation \eqref{Sys1_fun}, taken to hold in the sense of $L^2(0,T;V')$, can be shown to be equivalent to formulation \eqref{Sys1}.  In particular, the pressure gradient can be recovered using a corollary of a deep result of G. de Rham.  The corollary states that, for any distribution $g$, the equality $g=\nabla p$ holds for some distribution $p$ if and only if $\left<g, w\right> = 0$ for all $w \in \mathcal {V}$. See \cite{Wang_1993} for an elementary proof of the corollary.  

We recall the Agmon inequalities in 3D (see, e.g., \cite{Constantin_Foias_1988,Agmon}).  Namely, for any $\phi \in D(A)$, 
\begin{equation}
  \label{Agmon}
     \Vert\phi\Vert_{L^{\infty}} \leq C\Vert\nabla\phi\Vert_{L^2}^{1/2} \Vert A\phi\Vert_{L^2}^{1/2}
     \qquad\text{and}\qquad
     \Vert\phi\Vert_{L^{\infty}} \leq C\Vert\phi\Vert_{L^2}^{1/4} \Vert A\phi\Vert_{L^2}^{3/4}.
\end{equation}
 
The following Aubin-Lions Compactness Lemma is needed in order to construct solutions for \eqref{Sys1}.
\begin{Lemma}
\label{L2} Let $X$, $Y$, and $Z$ be separable, reflexive Banach spaces, where $X$ is compactly embedded in $Y$, and $Y$ is continuously embedded in $Z$.  
Let $T>0$, $p \in (1, \infty)$ and let $\{f_{n}(t, \cdot)\}_{n=1}^{\infty}$ 
be a bounded sequence in $L^{p}([0, T]; X)$
such that $\{\partial f_{n}/\partial t\}_{n=1}^{\infty}$ is bounded in $L^{p}([0, T]; Z)$. Then $\{f_{n}\}_{n=1}^{\infty}$ has a strongly convergent subsequence in $C([0, T]; Y)$. 
\end{Lemma} 

Typically in the theory of the Navier-Stokes equations, to write $\frac12\frac{d}{dt}\Vert u\Vert_{L^2}^2 = (u_t,u)$, one needs $u_t\in L^2(0,T;V')$, $u\in L^2(0,T;V)$ using the Lions-Magenes Lemma (cf., \cite[p. 176]{Temam_2001_Th_Num} or \cite[Corollary 7.3]{Robinson_2001}).  However, in our context, we have $u,\,u_t\in L^2(0,T;V)$.  Therefore, the following lemma is useful.  The proof is a straight-forward exercise via mollification in time, and follows closely the proofs of the Lions-Magenes Lemma in the aforementioned references.
\begin{Lemma}\label{Lemma_Lions_Magenes_type} 
Let $V$ be a Hilbert space with norm $\|\cdot\|$ and inner product $((\cdot,\cdot))$.  Suppose that for some $T>0$, $u\in L^2(0,T;V)$ and $u_t\in L^2(0,T;V)$.  Then the following equality holds in the scalar distribution sense on $(0,T)$.
\begin{align*}
 \frac12\frac{d}{dt}\|u\|^2 = ((u_t,u)).
\end{align*}
\end{Lemma} 

For the sake of completeness, we state the following uniform Gr\"onwall's inequality, proved in \cite{Jones_Titi_1992} (see also \cite{Farhat_Lunasin_Titi_2017_Horizontal} and the references therein), which will be used frequently throughout the paper. 
\begin{Lemma}
\label{L3}
Suppose that $Y(t)$ is a locally integrable and absolutely continuous function that satisfies the following: 
$$\frac{d Y}{d t} + \alpha(t) Y \leq \beta(t), \quad\text{ a.e. on } (0, \infty), $$
such that 
$$\liminf_{t \to \infty} \int_{t}^{t+\tau} \alpha(s)\,ds \geq \gamma, \quad\quad\quad \limsup_{t \to \infty} \int_{t}^{t+\tau} \alpha^{-}(s)\,ds < \infty, $$
and 
$$\lim_{t \to \infty} \int_{t}^{t+\tau} \beta^{+}(s)\,ds = 0, $$
for fixed $\tau > 0$, and $\gamma > 0$, 
where 
$\alpha^{-} = \max\{-\alpha, 0\}$ 
and 
$\beta^{+} = \max\{\beta, 0\}$. 
Then, $Y(t) \to 0$ at an exponential rate as $t \to \infty$.
\end{Lemma}

We record the following local well-posedness result for strong solutions to  equations \eqref{NSE}  (see, e.g., \cite{Temam_1995_Fun_Anal}, Theorem 4.2).
\begin{Theorem}\label{thm_NSE_local_reg}
 Let $\tilde u_0\in H^s\cap V$, $f\in L^2(0,T;H^{s-1}\cap H)$ for some $s\geq0$.  Then there exists a $T>0$ and unique solution $(\tilde u,\tilde p)$ to \eqref{NSE} such that $\tilde u\in C([0,T];H^s\cap V)\cap  L^2(0,T;H^{s+1}\cap V)$.
\end{Theorem}

\subsection{Main results of the paper}
\label{subsec2-2}
We first define the weak and strong solutions to system \eqref{Sys1}, respectively. 
\begin{Definition}
\label{Def1}
Let $T>0$ be arbitrary. Suppose $u_0 \in V$, $w_0 \in H$, and $f \in L^2(0, T; H)$. We call the pair $(u, w)$ a \textit{weak solution} on the time interval $[0, T]$ to system \eqref{Sys1}, if $u \in C(0, T; V)$, $u_t\in L^{2}(0, T; V)$, $w \in C_{w}(0, T; H) \cap L^2(0, T; V)$, $w_t\in L^{2}(0, T; H^{-1})$, and moreover, $(u, w)$ satisfies system \eqref{Sys1} in the weak sense, i.e., 
\begin{equation}
       \left\{
         \begin{aligned}
           &\alpha^2 ((u_{t}, \psi))
           +
           (u_{t}, \psi)
           +
           ((u, \psi))
           + 
           \left<w \times u, \psi\right> 
           = 
           \left(f, \psi\right),
           \\&
           \left<w_{t}, \psi\right> 
           + 
           ((w, \psi)) 
           - 
           \left<B(u, \psi), w\right>
           - 
           \left<\tilde{B}(w, u), \psi\right>
           =
           -\left(f, \nabla\times\psi\right),
         \end{aligned}
       \right.
    \label{Weak_def}
\end{equation}
holds for any $\psi \in L^2(0, T; V)$.  
\end{Definition}

Note that by taking $\psi = v\phi$ for $v\in V$ and $\phi\in C^1_c((0,T))$, it follows that formulation \eqref{Sys1_fun} is equivalent to formulation \eqref{Weak_def}, interpreted as an operator equation holding in an appropriate distributional sense.

\begin{Definition}
\label{Def2}
Let $T>0$ be an arbitrarily given time. Suppose $u_0 \in  V$, $w_0 \in V$, and $f \in L^2(0, T; H_{\text{curl}})$. We call the pair $(u, w)$ a \textit{strong solution} on the time interval $[0, T]$ to system \eqref{Sys1}, if it is a weak solution as in Definition~\ref{Def1} and satisfies additionally $w \in C([0, T]; V) \cap L^{2}(0, T; D(A))$, and $w_t\in L^{2}(0, T; H)$. 
\end{Definition}

The following theorem provides the global existence and uniqueness of weak solution to system \eqref{Sys1}.
\begin{Theorem}
\label{T1}
Suppose $u_0 \in V$, $w_0 \in H$, and $f \in L^2(0, T; H)$. Then, the velocity-vorticity-Voigt system \eqref{Sys1} possesses a unique global weak solution $(u, w)$ in the sense of Definition~\ref{Def1} that satisfies $\nabla\cdot w = 0$.  Moreover, the following energy equality holds.
\begin{align}\label{energy_equality_u}
 \alpha^2\|\nabla u(t)\|_{L^2}^2 
 + \|u(t)\|_{L^2}^2 
 +2\int_0^t\|\nabla u(s)\|_{L^2}^2\,ds
 =
 \alpha^2\|\nabla u_0\|_{L^2}^2 
 + \|u_0\|_{L^2}^2 
 +2\int_0^t(u(s),f(s))\,ds
\end{align}

\end{Theorem}

The next theorem is about the global existence and uniqueness of strong solution to system \eqref{Sys1}, 
as well as the higher-order regularity of the solution. 
\begin{Theorem}
\label{T2}
For the initial data $u_0 \in V$, $w_0 \in V$, and $f \in L^2(0, T; H_{\text{curl}})$,
there exists a unique strong solution $(u, w)$ in the sense of Definition~\ref{Def2}. 
Moreover, if we further assume that the initial data $u_0 \in H^s\cap V$, $w_0 \in H^s\cap V$, and $f \in L^2(0, T; H_{\text{curl}}^{s-1})$ for $s\geq2$, $s\in\mathbb{N}$, 
then, the solution $u \in C_{w}(0, T; H^s\cap V)$ and $w \in C_{w}(0, T; H^s\cap V) \cap L^2(0, T; H^{s+1}\cap V)$. 
\end{Theorem}

The following theorem relates the quantity $w$ in \eqref{Sys1_w} to the vorticity $\omega = \nabla\times u$. 
\begin{Theorem}
\label{T3}
Denote by $\omega:=\nabla\times u$ the vorticity of the flow and let $u_0 \in H^4\cap V$, $f\in H^{2}$. 
Then, we have 
$$\Vert\omega(t) - w(t)\Vert_{L^2}^2
+\alpha^2\Vert\nabla\omega(t) - \nabla w(t)\Vert_{L^2}^2 
+\int_0^t\Vert\nabla\omega(s) - \nabla w(s)\Vert_{L^2}^2\,ds
\leq 
C_0e^{Ct} + \frac{\tilde{K}\alpha^2}{C}(e^{Ct} - 1),$$
where $C_0$ depends on the initial data and $\tilde{K}$ is explained in the proof. 
If we further assume $w_0 = \nabla\times u_0$,  
then, $$\Vert \omega(t) - w(t)\Vert_{L^2}^2 
+ \alpha^2\Vert\nabla w(t) - \nabla\omega(t)\Vert_{L^2}^2 
+\int_0^t\Vert\nabla\omega(s) - \nabla w(s)\Vert_{L^2}^2\,ds
\leq 
K\alpha^2(e^{Ct} - 1),$$
for a.e. $t>0$, i.e., $\Vert w - \omega\Vert_{L^\infty(0,T;L^2)}\sim\mathcal{O}(\alpha)$ and $\Vert w - \omega\Vert_{L^2(0,T;V)}\sim\mathcal{O}(\alpha)$. 
In particular, we have $\Vert w - \omega\Vert_{L^\infty(0,T;L^2)} \to 0$ and 
$\Vert w - \omega\Vert_{L^2(0,T;V)} \to 0$ as $\alpha \to 0$.
\end{Theorem}

The next theorem describes the relation between the velocity-vorticity-Voigt equations \eqref{Sys1} and the 3D Navier-Stokes equations \eqref{NSE}.
\begin{Theorem}
\label{T4}
Denote by $\tilde{\omega}:=\nabla\times\tilde{u}$ the vorticity of $\tilde{u}$ in \eqref{NSE_u} and let $u_0$, $f$, and $T>0$ be the same as in Theorem~\ref{T3}, and set $w_0 = \nabla\times u_0$ and $\tilde{u}_0=u_0$. Then, for any $\alpha\in(0,1]$,
$$\Vert\omega(t) - \tilde{\omega}(t)\Vert_{L^2}^2 + \Vert u(t) - \tilde{u}(t)\Vert_{L^2}^2 + \alpha^2\Vert\nabla u(t) - \nabla\tilde{u}(t)\Vert_{L^2}^2 
+\int_0^t\Vert\nabla u(s) - \nabla\tilde{u}(s)\Vert_{L^2}^2\,ds
\leq 
C\alpha^2$$
for a.e. $t>0$ in the interval of existence of the solution to \eqref{NSE}, say, up to $T>0
$ and the constant $C$ depends on $\Vert\tilde{u}\Vert_{H^3}$, $\Vert u\Vert_{H^3}$, as well as $\Vert f\Vert_{H_{\text{curl}}^2}$. 
In particular, we have $\Vert \omega - \tilde{\omega}\Vert_{L^\infty(0,T;H)} \to 0$, $\Vert u - \tilde{u}\Vert_{L^\infty(0,T;H)} \to 0$, and $\Vert u - \tilde{u}\Vert_{L^2(0,T;V)} \to 0$ as $\alpha \to 0$. 
\end{Theorem}

\begin{Remark}
  \label{R1}
     We point out that the global well-posedness of \eqref{Sys1} still holds if we remove the divergence-free condition on the initial data $w_0$, i.e., we only assume $w_0\in L^2$ for the weak solution, and $w_0\in H^1$ for the strong solution. Consequently, we obtain the weak solutions $w\in C_{w}(0, T; L^2)\cap L^2(0, T; H^1)$ and the strong solution $w\in C_{w}(0, T; H^1)\cap L^2(0, T; H^2)$ by modifying the proofs in Section~\ref{sec3} and Section~\ref{sec4} accordingly. However, Theorem~\ref{T3} is no longer valid. Also, the assumptions that $\tilde{u}_0=u_0$ and $w_0 = \nabla\times u_0$ can be removed at the cost of obtaining ``convergence up to and error'' as $\alpha\rightarrow 0$.  For the sake of clarity, we do not include these details.
\end{Remark}

The above result yields the following blow-up criterion, which looks identical to the blow-up criterion for the 3D Euler-Voigt and 3D Navier-Stokes equations \cite{Larios_Titi_2009,Larios_Petersen_Titi_Wingate_2015} (see also \cite{Khouider_Titi_2008}).
\begin{Corollary}\label{blow_up}
 Assume the hypotheses and the notation of Theorem \ref{T4}.  Suppose that there is a $T>0$ and an $\epsilon>0$ such that
\begin{align}\label{blowup_criterion}
 \sup_{t\in[0,T]}\limsup_{\alpha\to0^+}\alpha\|\nabla u\|_{L^2}\geq\epsilon>0.
 \end{align}
Then solutions to the 3D Navier-Stokes equations with initial data $u_0$ develop a singularity on the time interval $[0,T]$, in the sense that there does not exist a strong solution $\tilde u\in L^2(0,T;V\cap H^2)\cap C([0,T],V)$.
\end{Corollary}

\section{Proof of Theorem~\ref{T1}}
\label{sec3}
In this section, we provide the construction of weak solution to system \eqref{Sys1} via Galerkin approximation. We also show that the obtained global weak solution is unique. 

\subsection{\em Proof of global existence of weak solutions}
\label{subsec3-1}
Consider the following finite-dimensional Galerkin ODE system for \eqref{Sys1}.
\begin{subequations}\label{ODE}
 \begin{empheq}[left=\empheqlbrace]{align}
   \label{ODEu}
     &\frac{d}{dt}(u_{M} + \alpha^2 Au_{M})
     +
     \nu Au_{M}
     +
     P_{M}P_\sigma(w_{M}\times u_{M})
     =
     f_{M}, 
     \\
     \label{ODEw}
     &\frac{d w_{M}}{dt}
     +
     A w_{M}
     +
     P_{M}B(u_{M}, w_{M})
     -
     P_{M}{\barB}(w_{M}, u_{M})
     =
     \nabla\times f_{M},
   \end{empheq}
\end{subequations}
with initial data $u_{M}(0) = P_{M}u_0$, $w_{M}(0) = P_{M}w_0$ and forcing $f_{M} = P_{M}f$. 

Notice that all the terms in \eqref{ODEu} and \eqref{ODEw} except the time-derivatives are at most quadratic, and thus, they are locally Lipschitz continuous. 
Therefore, by the Picard-Lindel\"of Theorem, we know that there exists a solution up to some time $T_{M} > 0$. 

Next we take inner-products with the above two equations by $u_{M}$ and $w_{M}$, respectively, integrate by parts, and obtain 
\begin{subequations}\label{E}
 \begin{empheq}[left=\empheqlbrace]{align}
   \label{Eu}
     &
     \frac{1}{2}\frac{d}{dt}
     \left( \Vert u_{M}\Vert_{L^2}^2 + \alpha^2\Vert\nabla u_{M}\Vert_{L^2}^2 \right)
     +
     \Vert\nabla u_{M}\Vert_{L^2}^2
     =
     \int_{\mathbb{T}^3}u_{M}\cdot f_{M}\,dx,  
     \\
     \label{Ew}
     &\frac{1}{2}\frac{d}{dt}\Vert w_{M}\Vert_{L^2}^2
     +
     \Vert\nabla w_{M}\Vert_{L^2}^2
     =
     \int_{\mathbb{T}^3}(w_{M}\cdot\nabla) u_{M}\cdot w_{M}\,dx
     +
     \int_{\mathbb{T}^3}(\nabla\times f_{M})\cdot w_{M}\,dx,
   \end{empheq}
\end{subequations}
where we used $\nabla\cdot u_{M} = 0$ in \eqref{Ew}. 
Then, by applying Cauchy-Schwarz inequality to the right side of \eqref{Eu}, we obtain
\begin{align}
     &
     \frac{d}{dt}
     \left( \Vert u_{M}\Vert_{L^2}^2 + \alpha^2\Vert\nabla u_{M}\Vert_{L^2}^2\right)
     +
     2\Vert\nabla u_{M}\Vert_{L^2}^2
     \leq
     2\Vert f_{M}\Vert_{L^2}\Vert u_{M}\Vert_{L^2}
     \leq
     \Vert f_{M}\Vert_{L^2}^2
     +
     \Vert\nabla u_{M}\Vert_{L^2}^2. 
     \label{Ineq1_u}
\end{align}
Thus, for any $T\in(0,T_M]$, integrating in time, it holds that for a.e. $t\in(0,T)$,
\begin{align*}
     \Vert u_{M}(t)\Vert_{L^2}^2 + \alpha^2\Vert\nabla u_{M}(t)\Vert_{L^2}^2 
     &\leq
     C\Vert f_M\Vert_{L^2(0,T,V)}^2
     +
      \Vert u_{M}(0)\Vert_{L^2}^2 + \alpha^2\Vert\nabla u_{M}(0)\Vert_{L^2}^2
     \\&\leq
     C\Vert f\Vert_{L^2(0,T;H)}^2
     +
      \Vert u_0\Vert_{L^2(0,T;H)}^2 + \alpha^2\Vert\nabla u_0\Vert_{L^2}^2
     =:K_{T}.
\end{align*}
Since the right-hand side is finite, $u_M$ can be extended beyond $T_M$, so that the above inequality holds for arbitrary $T>0$ and a.e. $t\in(0,T)$.  
 In particular, the interval of existence is independent of $M$.  Moreover, $\left\{u_M\right\}_{M=1^\infty}$ is uniformly bounded in $L^\infty(0,T;V)$. Using the Banach-Alaoglu Theorem, and extracting a subsequence if necessary (which we relabel if necessary and still denote by $u_M$), we obtain a $u\in L^{\infty}(0, T; V)$ such that $u_M$ converges to $u$ in the weak-$*$ sense of $L^{\infty}(0, T; V)$.

Next, we estimate the right side of \eqref{Ew} and obtain 
\begin{align}
\label{wM_L2_bound}
     \frac12\frac{d}{dt}\Vert w_{M}\Vert_{L^2}^2
     +
     \Vert\nabla w_{M}\Vert_{L^2}^2
     &\leq
     C\Vert\nabla u_{M}\Vert_{L^2} \Vert w_{M}\Vert_{L^2}^{1/2} \Vert\nabla w_{M}\Vert_{L^2}^{3/2}
     +
     \tfrac12\Vert\nabla\times f_{M}\Vert_{L^2}^2
     +
     \tfrac12\Vert w_{M}\Vert_{L^2}^2
     \notag
     \\&
     \leq
     C\Vert\nabla u_{M}\Vert_{L^2}^4 \Vert w_{M}\Vert_{L^2}^{2}
     +
     \tfrac12\Vert\nabla w_{M}\Vert_{L^2}^2
     +
     \tfrac12\Vert\nabla\times f_{M}\Vert_{L^2}^2
     +
     \tfrac12\Vert w_{M}\Vert_{L^2}^2
\end{align}
where we used Lemma~\ref{L1}. 
Then, after rearranging and using Gr\"onwall's inequality, it follows that
\begin{align*}
     \Vert w_{M}(t)\Vert_{L^2}^2
     &
     \leq
     \bar{K}_{T}
     := 
e^{TC_{K_{T}}}\Vert w_{M}(0)\Vert_{L^2}^2 + \int_0^Te^{C_{K_{T}}(t-s)}\Vert\nabla\times f(s)\Vert_{L^2}^2\,ds<\infty,
\end{align*}
from which we conclude that $w_{M}$ can be extended beyond $T_{M}$ up to any $t < T$ and is uniformly bounded in $L^{\infty}(0, T; L^2)$.  Using this fact and integrating \eqref{wM_L2_bound} on $[0,T]$, we also find that $w_{M}$ is uniformly bounded in $L^{2}(0, T; H^1)$. By similar arguments as above for $u_{M}$, we extract a weak-$*$ convergent subsequence, still denoted by $w_{M}$, with limit $w \in L^{\infty}(0, T; L^2)\cap L^{2}(0, T; H^1)$. 

Also, by taking the divergence of \eqref{ODEw} and denoting $v_{M} := \nabla\cdot w_{M}$, we obtain
\begin{align}
     &
     \frac{d v_{M}}{d t}
     +
     A v_{M}
     +
     P_{M}B(u_{M}, v_{M})
     =
     0. 
  \label{Div-w}
\end{align}
Multiplying the above equation by $v_{M}$ and using \eqref{embd2} and \eqref{Poincare}, we obtain
\begin{align*}
     &
     \frac{1}{2}\frac{d}{d t} \Vert v_{M}\Vert_{L^2}^2
     +
     \lambda_1\Vert v_{M}\Vert_{L^2}^2
     \leq
     0. 
\end{align*}
Then, Gr\"onwall's inequality the implies for a.e. $t\in(0,T)$, 
$$\Vert v_{M}(t)\Vert_{L^2} \leq \Vert v_{M}(0)\Vert_{L^2}e^{-\lambda_1t}. $$ 
If the initial data $w_0$ is in $H$, we have $v_{M}(0) = 0$, 
which implies $\nabla\cdot w_{M} = v_{M} = 0$ in $L^2([0,T];L^2)$.  Since $L^2([0,T];H)$ is closed in $L^2([0,T];L^2)$, this also implies $\nabla\cdot w=0$, so long as $\nabla\cdot w_0=0$. 
Namely, we have $w\in L^{\infty}(0, T; H)\cap L^{2}(0, T; V)$. 
Note that \eqref{Div-w} also implies 
\begin{align*}
     &
     \Vert v_{M}(T)\Vert_{L^2}^2
     +
     2\int_{0}^{T}\Vert\nabla v_{M}\Vert_{L^2}^2\,dt
     \leq
     \Vert v_{M}(0)\Vert_{L^2}^2,  
\end{align*}
from which we obtain that $v_{M}$ is uniformly bounded in $L^{2}(0, T; H^1)$. 

We consider the pair $(u, w)$ as our candidate solution. 
Next, we obtain bounds on 
$d u_{M}/dt$ in $L^{2}(0, T; V)$ and bounds on $d w_{M}/dt$ in $L^{2}(0, T; V')$, uniformly with respect to $M$.
Note that
\begin{subequations}\label{Time}
 \begin{empheq}[left=\empheqlbrace]{align}
   \label{Time_u}
     &
     (I + \alpha^2A)\frac{d u_{M}}{d t}
     =
     -
     Au_{M}
     -
     P_{M}(u_{M}\times w_{M})
     +
     f_{M},
     \\&
   \label{Time_w}
     \frac{d w_{M}}{dt}
     =
     -A w_{M}
     -
     P_{M}B(u_{M}, w_{M})
     +
     P_{M}{\barB}(w_{M}, u_{M})
     +
     \nabla\times f_{M}.
    \end{empheq}
\end{subequations}
Note in \eqref{Time_u} that $-Au_{M}$ is uniformly bounded in $L^{2}(0, T; V')$ 
due to the fact that $u_{M}$ is also uniformly bounded in $L^{2}(0, T; V)$. 
On the other hand, by Lemma~\ref{L1}, we have 
\begin{align*}
\left|\int_{\mathbb{T}^3}P_M(u_{M}\times w_{M})\cdot \phi\,dx\right|
&=
     \left|\int_{\mathbb{T}^3}(u_{M}\times w_{M})\cdot P_M\phi\,dx\right|
     \leq
     C\Vert u_{M}\Vert_{L^2}^{1/2} \Vert\nabla u_{M}\Vert_{L^2}^{1/2} \Vert w_{M}\Vert_{L^2} \Vert P_M\phi\Vert_{H^1}
     \\&\leq
     C\Vert u_{M}\Vert_{L^2}^{1/2} \Vert\nabla u_{M}\Vert_{L^2}^{1/2} \Vert w_{M}\Vert_{L^2} \Vert \phi\Vert_{H^1},
\end{align*}
for all test functions $\phi\in V$.
Thus, $P_{M}(u_{M}\times w_{M})$ is also uniformly bounded in $L^{2}(0, T; V')$. 
It is easily seen that, 
uniformly in $M$, $f_{M}$ and $\nabla\times f_{M}$ are bounded in $L^{2}(0, T; V')$, 
as well as $Au_{M}$, $A w_{M}$.
Therefore, 
$$(I + \alpha^2 A)\frac{d u_{M}}{d t} \text{\quad is uniformly bounded in\quad} L^2(0, T; V').$$
By inverting the Helmholtz operator $(I+\alpha^2 A)$ with respect to zero-mean, periodic boundary conditions, we obtain 
$$\frac{d u_{M}}{d t} \text{\quad is uniformly bounded in\quad} L^2(0, T; V).
$$

Next, notice in \eqref{Time_w} that $\Delta w_{M}$ is bounded in $L^{2}(0, T; H^{-1}) \subset L^{2}(0, T; V')$, 
and for the nonlinear terms, 
we integrate by parts and use Lemma~\ref{L1} in order to get 
\begin{align*}
     \left|\int_{\mathbb{T}^3}P_M((u_{M}\cdot\nabla) w_{M})\cdot \phi\,dx\right|
     &=
     \left|\int_{\mathbb{T}^3}(u_{M}\cdot\nabla) w_{M}\cdot P_M\phi\,dx\right|
     \leq
     \Vert u_{M}\Vert_{L^2}^{1/2} \Vert\nabla u_{M}\Vert_{L^2}^{1/2} \Vert w_{M}\Vert_{H^1} \Vert\phi\Vert_{H^1}
\end{align*}
for all test functions $\phi\in V$. 
Therefore,  $-P_{M}B(u_{M}, w_{M})$ is uniformly bounded in $L^{2}(0, T; V')$. 
Similar estimates show that $P_{M}{\barB}(w_{M}, u_{M})$ is also bounded in $L^{2}(0, T; V')$, uniformly in $M$. 
Thus, we get 
$$\frac{d w_{M}}{d t} \text{\quad is uniformly bounded in\quad} L^2(0, T; V').$$
Also, by the bounds we obtained above and Lemma~\ref{L2}, 
there is a subsequence, still labeled as $(u_{M}, w_{M})$ that satisfies
\begin{subequations}
\begin{align}
     \label{Convergence1}
     u_{M} \rightarrow u \text{\quad strongly in \quad} L^2(0, T; V) \text{\quad and\quad} &w_{M} \rightarrow w \text{\quad strongly in \quad} L^2(0, T; H),
     \\
     \label{Convergence2}
     &w_{M} \rightharpoonup w \text{\quad weakly in\quad} L^2(0, T; V),
     \\
     \label{Convergence3}
     u_{M} \rightharpoonup u \text{\quad weak-$\ast$ in\quad} L^{\infty}([0, T]; V) \text{\quad and\quad} &w_{M} \rightharpoonup w \text{\quad weak-$\ast$ in\quad} L^{\infty}(0, T; H),
\end{align}
\end{subequations}
for all $T>0$.
Hence, for $0<\bar{T}\leq\tilde{T}<T$, by taking inner products of \eqref{ODEu} and \eqref{ODEw} with test function $\psi \in C_{c}^1([0, \tilde{T}); V)$, integrating in time over $[0, \bar{T}]$, and integrating by parts, we obtain
\begin{subequations}\label{Weak}
 \begin{empheq}[left=\empheqlbrace]{align}
   \label{Weak_u}
     &-\int_{0}^{\bar{T}}(u_{M}, \psi_{t})\,dt
     +
     (u_{M}(\cdot, \bar{T}), \psi(\cdot, \bar{T}))
     -
     (u_{M}(\cdot, 0), \psi(\cdot, 0)) \nonumber
     \\&\quad
     -
     \alpha^2\int_{0}^{\bar{T}}((u_{M}, \psi_{t}))\,dt
     -
     \alpha^2((u_{M}(\cdot, \bar{T}), \psi(\cdot, \bar{T})))
     +
     \alpha^2((u_{M}(\cdot, 0), \psi(\cdot, 0))) \nonumber
     \\&\quad\quad
     +
     \int_{0}^{\bar{T}}((u_{M}, \psi))\,dt
     +
     \int_{0}^{\bar{T}}\left<(w_{M}\times u_{M}), P_{M}\psi\right>\,dt
     =
     \int_{0}^{\bar{T}}(f_{M}, \psi)\,dt, 
     \\
     \label{Weak_w}
     &-\int_{0}^{\bar{T}}(w_{M}, \psi_{t})\,dt
     +
     (w_{M}(\cdot, \bar{T}), \psi(\cdot, \bar{T}))
     -
     (w_{M}(\cdot, 0), \psi(\cdot, 0))
     +
     \int_{0}^{\bar{T}}((u_{M}, \psi))\,dt \nonumber
     \\&\quad
     +
     \int_{0}^{\bar{T}}\left<B(u_{M}, w_{M}), P_{M}\psi\right>\,dt
     -
     \int_{0}^{\bar{T}}\left<{\barB}(w_{M}, u_{M}), P_{M}\psi\right>\,dt
     =
     \int_{0}^{\bar{T}}(f_{M}, \nabla\times\psi)\,dt. 
   \end{empheq}
\end{subequations}
Using the standard arguments from the theory of the Navier-Stokes equations (see, e.g., \cite{Constantin_Foias_1988, Temam_2001_Th_Num}), 
we have that each of the integrals in \eqref{Weak_u} and \eqref{Weak_w} converges to the time integral of the corresponding term in \eqref{Weak}. 
For the sake of completeness, we provide the details below. 
First, convergence of integrals of the linear terms follows from \eqref{Convergence1} as well as $f\in H$. 
The choice of $\psi$ and \eqref{Convergence1} imply the convergence of the boundary terms at $t=0$ and $t=\bar{T}$ in both \eqref{Weak_u} and \eqref{Weak_w}. 
Regarding the nonlinear term in \eqref{Weak_u}, we have 
\begin{align*}
     &\left| \int_{0}^{\bar{T}}\left<(w_{M}\times u_{M}), P_{M}\psi\right>\,dt - \int_{0}^{\bar{T}}\left<(w\times u), \psi\right>\,dt \right|
     \\&
     \leq
     \int_{0}^{\bar{T}}\left|\left<(w_{M}\times (u_{M}-u)), P_{M}\psi\right>\right|\,dt
     +
     \int_{0}^{\bar{T}}\left|\left<(w_{M}-w)\times u, P_{M}\psi\right>\right|\,dt
     \\&\quad
     +
     \int_{0}^{\bar{T}}\left|\left<w\times u, P_{M}\psi - \psi\right>\right|\,dt
     \\&
     \leq
     \Vert w_{M}\Vert_{L^2(L^2)} \Vert u-u_{M}\Vert_{L^2(H^1)} \Vert\psi\Vert_{L^{\infty}(L^2)}^{1/2} \Vert\psi\Vert_{L^{\infty}(H^1)}^{1/2}
     \\&\quad
     +
     \Vert w_{M}-w\Vert_{L^2(L^2)} \Vert u\Vert_{L^2(H^1)} \Vert\psi\Vert_{L^{\infty}(L^2)}^{1/2} \Vert\psi\Vert_{L^{\infty}(H^1)}^{1/2}
     \\&\quad
     +
     \Vert w\Vert_{L^2(H^1)} \Vert u\Vert_{L^{\infty}(L^2)}^{1/2}\Vert u\Vert_{L^{\infty}(H^1)}^{1/2} \Vert P_{M}\psi - \psi\Vert_{L^2(L^2)}
\end{align*}
which converges to $0$ in view of \eqref{Convergence1} and the uniform bounds on $w_{M}$ and $u$. 
As for the first nonlinear term in \eqref{Weak_w}, we use \eqref{embd1} and estimate as 
\begin{align*}
     &\left| \int_{0}^{\bar{T}}\left<B(u_{M}, P_{M}\psi), w_{M}\right>\,dt - \int_{0}^{\bar{T}}\left<B(u, \psi), w\right>\,dt\right|
     \\&
     \leq
     \int_{0}^{\bar{T}}\left|\left<B(u_{M}, P_{M}\psi - \psi), w_{M}\right>\right|\,dt
     +
     \int_{0}^{\bar{T}}\left|\left<B(u_{M}-u, \psi), w_{M}\right>\right|\,dt
     \\&\quad
     +     
     \left|\int_{0}^{\bar{T}}\left<B(u_{M}, \psi), w_{M}-w\right>\,dt\right|
     \\&
     \leq
     \Vert u_{M}\Vert_{L^2(H^1)} \Vert P_{M}\psi-\psi\Vert_{L^{\infty}(H^1)} \Vert w_{M}\Vert_{L^2(L^2)}^{1/2} \Vert w_{M}\Vert_{L^2(H^1)}^{1/2}
     \\&\quad
     +
     \Vert u_{M}-u\Vert_{L^2(H^1)} \Vert\psi\Vert_{L^{\infty}(H^1)} \Vert w_{M}\Vert_{L^2(L^2)}^{1/2} \Vert w_{M}\Vert_{L^2(H^1)}^{1/2}
     \\&\quad
     +
     \Vert u_{M}\Vert_{L^{2}(H^1)} \Vert \psi\Vert_{L^{\infty}(H^1)} \Vert w_{M}-w\Vert_{L^2(L^2)}^{1/2} \Vert w_{M}-w\Vert_{L^2(H^1)}^{1/2}
     \\&
     \rightarrow 0 \text{\quad as\quad } M\rightarrow\infty,
\end{align*}
where we used \eqref{Convergence2}, the uniform boundedness of $w_{M}$ in $H$ and integrated by parts in the last integral. 
Finally, convergence of the second nonlinear term in \eqref{Weak_w} is obtained as 
\begin{align*}
      &\left| \int_{0}^{\bar{T}}\left<{\barB}(w_{M}, u_{M}), P_{M}\psi\right>\,dt - \int_{0}^{\bar{T}}\left<{\barB}(w, u), \psi\right>\,dt\right|
     \\&
     \leq
     \int_{0}^{\bar{T}}\left|\left<{\barB}(w_{M}, u_{M}), P_{M}\psi-\psi\right>\right|\,dt
     +
     \int_{0}^{\bar{T}}\left|\left<{\barB}(w_{M}, \psi), u_{M}-u\right>\right|\,dt
     \\&\quad
     +     
     \int_{0}^{\bar{T}}\left|\left<{\barB}(w_{M}-w, u), \psi\right>\right|\,dt
     \\&
     \leq
     \Vert w_{M}\Vert_{L^2(L^{2})}^{1/2}\Vert w_{M}\Vert_{L^2(H^{1})}^{1/2} \Vert u_{M}\Vert_{L^{2}(H^2)} \Vert P_{M}\psi-\psi\Vert_{L^{\infty}(L^2)}
     \\&\quad
     +
     \Vert w_{M}\Vert_{L^2(L^2)}^{1/2}\Vert w_{M}\Vert_{L^2(H^1)}^{1/2} \Vert\psi\Vert_{L^{\infty}(H^1)} \Vert u_{M}-u\Vert_{L^{2}(H^1)} 
     \\&\quad
     +
     \Vert w_{M}-w\Vert_{L^2(L^{2})} \Vert u\Vert_{L^{2}(H^2)} \Vert\psi\Vert_{L^{\infty}(L^2)}^{1/2} \Vert\psi\Vert_{L^{\infty}(H^1)}^{1/2}
     \\&
     \rightarrow 0 \text{\quad as\quad } M\rightarrow\infty,
\end{align*}
due to \eqref{Convergence2} and the uniform boundedness of $w_{M}$. 
Similar arguments also apply to equation \eqref{Div-w} for $v_{M}:=\nabla\cdot w_{M}$.
Namely, each term in 
\begin{align*}
     &-\int_{0}^{\bar{T}}(v_{M}, \psi_{t})\,dt
     +
     (v_{M}(\cdot, \bar{T}), \psi(\cdot, \bar{T}))
     -
     (v_{M}(\cdot, 0), \psi(\cdot, 0))
     \\&\quad
     +
     \int_{0}^{\bar{T}}((v_{M}, \psi))\,dt
     +
     \int_{0}^{\bar{T}}\left<B(u_{M}, v_{M}), P_{M}\psi\right>\,dt
     =
     0
\end{align*}
converges to the time integral of the following weak formulation of $v$ in view of \eqref{Div-w}, 
\begin{align*}
     &
     \left<v_{t}, \psi\right>
     +
     ((v, \psi))
     +
     \left<B(u, v), \psi\right>
     =
     0. 
\end{align*} 
Specifically, the linear terms converge due to the fact that $v_{M}$ is bounded in $L^{2}(0, T; H^1)$ from \eqref{Div-w}, 
while the convergence of the nonlinear term follows similar to those in \eqref{Weak_w}, i.e., 
after integration by parts, we have 
\begin{align*}
     &\left| \int_{0}^{\bar{T}}\left<B(u_{M}, P_{M}\psi), v_{M}\right>\,dt - \int_{0}^{\bar{T}}\left<B(u, \psi), v\right>\,dt\right|
     \\&
     \leq
     \int_{0}^{\bar{T}}\left|\left<B(u_{M}, P_{M}\psi - \psi), v_{M}\right>\right|\,dt
     +
     \int_{0}^{\bar{T}}\left|\left<B(u_{M}-u, \psi), v_{M}\right>\right|\,dt
     \\&\quad
     +     
     \left|\int_{0}^{\bar{T}}\left<B(u_{M}, \psi), v_{M}-v\right>\,dt\right|
     \\&
     \leq
     \Vert u_{M}\Vert_{L^2(H^1)} \Vert P_{M}\psi-\psi\Vert_{L^{\infty}(H^1)} \Vert v_{M}\Vert_{L^2(L^2)}^{1/2} \Vert v_{M}\Vert_{L^2(H^1)}^{1/2}
     \\&\quad
     +
     \Vert u_{M}-u\Vert_{L^2(H^1)} \Vert\psi\Vert_{L^{\infty}(H^1)} \Vert v_{M}\Vert_{L^2(L^2)}^{1/2} \Vert v_{M}\Vert_{L^2(H^1)}^{1/2}
     \\&\quad
     +
     \Vert u_{M}\Vert_{L^{2}(H^1)} \Vert \psi\Vert_{L^{\infty}(H^1)} \Vert v_{M}-v\Vert_{L^2(L^2)}^{1/2} \Vert v_{M}-v\Vert_{L^2(H^1)}^{1/2}
     \\&
     \rightarrow 0 \text{\quad as\quad } M\rightarrow\infty,
\end{align*}
where we also used $v_M\rightarrow v$ in $L^2(0, T; L^2)$ as $M\rightarrow \infty$, as well as \eqref{Convergence1}. 
Note that by the fact $v_{M}=0$ in $L^2([0,T];L^2)$, 
we have $v(t) = \nabla\cdot w(t) = 0$ for a.e. $t<T$. 
Now, all the above convergence is valid if we take $\psi=v\phi(t)$ where $v\in C^{\infty}$ and $\phi\in C^{1}({0, \tilde{T}})$. In particular, the convergence is valid for all $\phi\in \mathcal{D}([0,\tilde{T}])$, thus, \eqref{Weak_def} holds in the sense of distributions, which in turn implies that  \eqref{Sys1_fun} is valid as an equation of operators. 
Namely, in view of the embedding $V \hookrightarrow H \hookrightarrow V'$, we conclude that equations \eqref{Sys1} hold in the weak sense by Lemma~\ref{L2}, 
while the pressure term $p$ is recovered by the approach mentioned in Section~\ref{sec2}. 
Finally, by integrating in time over $[\bar{t}, \bar{T}]$ for $0\leq \bar{t} < \bar{T}$ and sending $\bar{T} \to \bar{t}$ one can use similar convergence arguments as above (c.f. \cite{Constantin_Foias_1988, Temam_2001_Th_Num}) and show that $(u, w)$ is in fact weakly continuous with respect to time in $V\times H$, so that the initial condition is satisfied in the weak sense. 
Since $T>0$ is arbitrary, the existence in Theorem~\ref{T1} is thus proven. 

\subsection{\em Proof of uniqueness}
\label{subsec3-2}
Suppose there exist two pairs of solution $(u, w, p)$ and $(\bar{u}, \bar{w}, \bar{p})$ to system \eqref{Sys1}, with the same initial data $u_0\in V$, $w_0\in H$ and the same forcing $f$ on their common time interval of existence $(0, T)$. 
By subtracting the equations of the two pairs of solution 
and denoting $\tilde{u} = u - \bar{u}$, $\tilde{w} = w - \bar{w}$, and $\tilde{p} =  p- \bar{p}$, 
we obtain 
\begin{equation}
       \left\{
         \begin{aligned}
           &(I - \alpha^2\Delta)\frac{\partial \tilde{u}}{\partial t} 
           -
           \Delta \tilde{u}
           + 
           \tilde{w} \times u 
           +
           \bar{w} \times \tilde{u}
           + 
           \nabla \tilde{p} 
           = 
           0,
           \\&
           \frac{\partial \tilde{w}}{\partial t} 
           - 
           \Delta \tilde{w} 
           + 
           (u\cdot\nabla) \tilde{w}
           +
           (\tilde{u}\cdot\nabla) \bar{w}
           - 
           (w\cdot\nabla) \tilde{u}
           -
           (\tilde{w}\cdot\nabla) \bar{u}
           =
           0,
           \\&
           \nabla \cdot \tilde{u} = 0,
           \\&
           \tilde{u}(\cdot, 0) = 0
           \quad\text{ and }\quad
           \tilde{w}(\cdot, 0) = 0. 
         \end{aligned}
       \right.
    \label{Diff}
\end{equation}
Multiplying the two equations by $\tilde{u}$ and $\tilde{w}$, respectively, 
integrating by parts over $\mathbb{T}^3$, and adding, we obtain 
\begin{align}
     &
     \frac{1}{2}\frac{d}{d t}\left(\Vert\tilde{u}\Vert_{L^2}^2 + \alpha^2\Vert\nabla\tilde{u}\Vert_{L^2}^2 + \Vert\tilde{w}\Vert_{L^2}^2\right)
     +
     \Vert\nabla\tilde{u}\Vert_{L^2}^2
     +
     \Vert\nabla\tilde{w}\Vert_{L^2}^2 \nonumber
     \\&
     =
     -
     \int_{\mathbb{T}^3} (\tilde{w}\times u)\cdot \tilde{u}\,dx
     -
     \int_{\mathbb{T}^3} (\tilde{u}\cdot\nabla)\bar{w}\cdot \tilde{w}\,dx
     -
     \int_{\mathbb{T}^3} (w\cdot\nabla)\tilde{u}\cdot \tilde{w}\,dx
     +
     \int_{\mathbb{T}^3} (\tilde{w}\cdot\nabla)\bar{u}\cdot \tilde{w}\,dx, 
   \label{Diff_eq}
\end{align}
where we used $\nabla\cdot u = \nabla\cdot \bar{u} = \nabla\cdot \tilde{u} = 0$. 
Next, we estimate the four terms on the right side of \eqref{Diff_eq}. 
By Lemma~\ref{L1}, H\"older's and Young's inequalities, 
the first integral is bounded by 
\begin{align*}
     \int_{\mathbb{T}^3} |u| |\tilde{u}| |\tilde{w}|\,dx
     &
     \leq
     C\Vert u\Vert_{L^2} \Vert\tilde{u}\Vert_{L^2}^{1/2} \Vert\nabla\tilde{u}\Vert_{L^2}^{1/2} \Vert\nabla\tilde{w}\Vert_{L^2}
     \leq
     C\Vert u\Vert_{L^2}^2 \Vert\tilde{u}\Vert_{L^2} \Vert\nabla\tilde{u}\Vert_{L^2}
     +
     \frac{1}{8}\Vert\nabla\tilde{w}\Vert_{L^2}^2
     \\&
     \leq
     C_{\alpha}\Vert u\Vert_{L^2}^2\left( \Vert\tilde{u}\Vert_{L^2}^2 + \alpha^2\Vert\nabla\tilde{u}\Vert_{L^2}^2 \right)
     +
     \frac{1}{8}\Vert\nabla\tilde{w}\Vert_{L^2}^2.
\end{align*}
By similar estimates, the second integral is bounded by
\begin{align*}
     \int_{\mathbb{T}^3} |\bar{w}| |\tilde{u}| |\nabla\tilde{w}|\,dx
     &
     \leq
     C\Vert\bar{w}\Vert_{H^1} \Vert\nabla\tilde{u}\Vert_{L^2} \Vert\nabla\tilde{w}\Vert_{L^2}
     \leq
     C_{\alpha} \Vert\bar{w}\Vert_{H^1}^2 \left(\alpha^2 \Vert\nabla\tilde{u}\Vert_{L^2}^2\right)
     +
     \frac{1}{8} \Vert\nabla\tilde{w}\Vert_{L^2}^2, 
\end{align*}
where we integrated by parts and used $\nabla\cdot\tilde{u} = 0$. 
By Lemma~\ref{L1}, the third integral is bounded by 
\begin{align*}
     \int_{\mathbb{T}^3} |w| |\nabla\tilde{u}| |\tilde{w}|\,dx
     &
     \leq
     C\Vert w\Vert_{H^1} \Vert\nabla\tilde{u}\Vert_{L^2} \Vert\tilde{w}\Vert_{L^2}^{1/2} \Vert\nabla\tilde{w}\Vert_{L^2}^{1/2}
     \leq
     C\Vert w\Vert_{H^1}^{4/3} \Vert\nabla\tilde{u}\Vert_{L^2}^{4/3} \Vert\tilde{w}\Vert_{L^2}^{2/3}
     +
     \frac{1}{8} \Vert\nabla\tilde{w}\Vert_{L^2}^2 
     \\&
     \leq
     C_{\alpha} \Vert w\Vert_{H^1}^2 \left(\alpha^2 \Vert\nabla\tilde{u}\Vert_{L^2}^2\right)
     +
     C_{\alpha} \Vert\tilde{w}\Vert_{L^2}^2
     +
     \frac{1}{8} \Vert\nabla\tilde{w}\Vert_{L^2}^2,  
\end{align*}
where we applied H\"older's and Young's inequalities. 
As for the last integral, we obtain the upper bound analogously as 
\begin{align*}
     \int_{\mathbb{T}^3} |\nabla\bar{u}| |\tilde{w}|^2\,dx
     &
     \leq
     C\Vert\bar{u}\Vert_{H^1} \Vert\tilde{w}\Vert_{L^2}^{1/2} \Vert\nabla\tilde{w}\Vert_{L^2}^{3/2}
     \leq
     C\Vert\bar{u}\Vert_{H^1}^4 \Vert\tilde{w}\Vert_{L^2}^{2}
     +
     \frac{1}{8} \Vert\nabla\tilde{w}\Vert_{L^2}^2.
\end{align*}
Summing up all the above estimates, 
we obtain 
\begin{align*}
     &
     \frac{d}{d t}\left(\Vert\tilde{u}\Vert_{L^2}^2 + \alpha^2\Vert\nabla\tilde{u}\Vert_{L^2}^2 + \Vert\tilde{w}\Vert_{L^2}^2\right)
     +
     \Vert\nabla\tilde{u}\Vert_{L^2}^2
     +
     \Vert\nabla\tilde{w}\Vert_{L^2}^2 
     \leq
     M_{\alpha}\left(\Vert\tilde{u}\Vert_{L^2}^2 + \alpha^2\Vert\nabla\tilde{u}\Vert_{L^2}^2 + \Vert\tilde{w}\Vert_{L^2}^2\right). 
\end{align*}
Here $M_{\alpha}:=M_{1} + M_{2}(\Vert w\Vert_{H^1}^2+\Vert\bar{w}\Vert_{H^1}^2)$, is such that $M_1$ depends on $\Vert\bar{u}\Vert_{H^1}^4$ and $\Vert u\Vert_{L^2}^2$, which are bounded, while $M_2$ is an absolute constant. Therefore, by Gr\"onwall's inequality and $w, \bar{w}\in L^{2}(0, T; V)$, we conclude that 
$$\Vert\tilde{u}(t)\Vert_{L^2}^2 + \Vert\tilde{w}(t)\Vert_{L^2}^2 = 0,$$ 
since $\tilde{u}_0 = \tilde{w}_0 = 0$. 
Namely, we have $u(t) = \bar{u}(t)$ and $w(t) = \bar{w}(t)$. Finally, setting $\psi = u\in L^2(0,T;V)$ in Definition \ref{Def1}, using Lemma \eqref{Lemma_Lions_Magenes_type}, and integrating in time, we obtain \eqref{energy_equality_u}.  
The proof of Theorem~\ref{T1} is thus complete.

\section{Proof of Theorem~\ref{T2}}
\label{sec4}
In this section, we show that system \eqref{Sys1} has a unique global strong solution and provide the {\it{a priori}} estimates for the higher order regularity of such solution $(u, w)$ to \eqref{Sys1}, with $u_0, w_0 \in V$. In view of Definition~\ref{Def2}, it suffices to prove the uniform boundedness of $w_{M}$ in $V$. 

\subsection{\em Proof of global existence of strong solutions}
\label{subsec4-1}

To begin, we multiply \eqref{ODEw} by $Aw_{M}$, respectively, integrate by parts over $\mathbb{T}^3$ and obtain 
\begin{align}
     \frac{1}{2}\frac{d}{d t}\Vert\nabla w_{M}\Vert_{L^2}^2
     +
     \Vert Aw_{M}\Vert_{L^2}^2 \nonumber
     &=
     \int_{\mathbb{T}^3} (u_{M}\cdot\nabla)w_{M}\cdot Aw_{M}\,dx
     -
     \int_{\mathbb{T}^3} (w_{M}\cdot\nabla)u_{M}\cdot Aw_{M}\,dx \nonumber
     \\&\quad
     -
     \int_{\mathbb{T}^3} (\nabla\times f_{M})\cdot Aw_{M}\,dx. 
   \label{Eq3}
\end{align}
Then, we estimate the three terms on the right side of \eqref{Eq3}. 
For the first term, we integrate by parts and apply Lemma~\ref{L1} as 
\begin{align*}
     \int_{\mathbb{T}^3} (u_{M}\cdot\nabla)w_{M}\cdot Aw_{M}\,dx
     &\leq
     C\Vert\nabla u_{M}\Vert_{L^2} \Vert\nabla w_{M}\Vert_{L^2}^{1/2} \Vert Aw_{M}\Vert_{L^2}^{3/2}
     \\&
     \leq
     C\Vert\nabla u_{M}\Vert_{L^2}^{4} \Vert\nabla w_{M}\Vert_{L^2}^{2}
     +
     \frac{1}{8}\Vert Aw_{M}\Vert_{L^2}^2. 
\end{align*}

By \eqref{Agmon}, the second term is bounded by 
\begin{align*}
     \int_{\mathbb{T}^3}|w_{M}| |\nabla u_{M}| |Aw_{M}|\,dx
     &\leq
     C\Vert w_{M}\Vert_{L^{\infty}} \Vert\nabla u_{M}\Vert_{L^2} \Vert Aw_{M}\Vert_{L^2}
     \leq     
     C \Vert\nabla u_{M}\Vert_{L^2} \Vert\nabla w_{M}\Vert_{L^2}^{1/2} \Vert Aw_{M}\Vert_{L^2}^{3/2} 
     \\&
     \leq
     C\Vert\nabla u_{M}\Vert_{L^2}^4\Vert\nabla w_{M}\Vert_{L^2}^2 
     +
     \frac{1}{8}\Vert Aw_{M}\Vert_{L^2}^2. 
\end{align*} 
By H\"older's inequality, the last term is bounded by 
\begin{align*}
     \int_{\mathbb{T}^3}|\nabla\times f_{M}| |Aw_{M}|\,dx
     &\leq
     \Vert\nabla\times f_{M}\Vert_{L^2} \Vert Aw_{M}\Vert_{L^2}
     \leq
     C\Vert\nabla\times f\Vert_{L^2}^2
     +
     \frac{1}{8}\Vert Aw_{M}\Vert_{L^2}^2. 
\end{align*}
Combining all the above estimates and denoting 
$$X_{M}(t) = \Vert\nabla w_{M}(t)\Vert_{L^2}^2 \text{\quad and \quad} Y_{M}(t) = \Vert Aw_{M}(t)\Vert_{L^2}^2,$$
for $0\leq t\leq T$, 
we get 
\begin{align*}
     \frac{d}{d t}X_{M}(t) + Y_{M}(t)
     \leq
     CX_{M}(t) + C\Vert f\Vert_{H_{\text{curl}}}^2. 
\end{align*}
Thus, the uniform bound of $u$ in $V$ and Gr\"onwall's inequality imply that $w_{M}$ is uniformly bounded in $L^{\infty}(0, T; V)$. 
Integrating in time over $[0, T]$, we also have that $w_{M}$ is uniformly bounded in $L^2(0, T; D(A))$. 
As for the time derivative of $w_{M}$ in \eqref{Time_w}, 
we use the fact that $w_{M}$ are bounded in $L^{2}(0, T; D(A))$ 
and all the nonlinear terms are bounded in $L^{2}(0, T; L^2)$, and conclude that  
$$\frac{d w_{M}}{d t} \text{\quad is uniformly bounded in\quad} L^2(0, T; H).$$ 
A standard simple argument (see, e.g., \cite{Temam_2001_Th_Num}) then shows that $dw/dt\in L^2(0, T; H)$, thanks to the convergence properties of $w_M$ proven above.
Therefore, we obtain the global existence of strong solutions. 

\subsection{\em Proof of higher regularity}
\label{subsec4-3}
In this subsection, we first provide the $H^2$ and $H^3$ {\it{a priori}} estimates for $u$ and $w$, 
then we obtain the $H^s$ bounds on $u$ and $w$ for $s\geq 4$. 
\subsubsection{\em $H^2$ bounds}
\label{subsubsec4-3-1}
First we provide the $H^2$ {\it a priori} estimates for $u$ with initial data $u_0\in D(A)$. We work formally for the sake of clarity, we point out that the following estimates can be justified at the Galerkin level following similar arguments to those in Section~\ref{subsec3-1} and Section~\ref{subsec4-1}. 
We begin by multiplying \eqref{Sys1_u} by $-\Delta u$, integrating by parts over $\mathbb{T}^3$, and obtain
\begin{align}
     &
     \frac{1}{2}\frac{d}{d t}\left(\Vert\nabla u\Vert_{L^2}^2 + \alpha^2\Vert \Delta u\Vert_{L^2}^2\right)
     +
     \Vert \Delta u\Vert_{L^2}^2 
     =
     \int_{\mathbb{T}^3} w\times u\cdot\Delta u\,dx
     +
     \int_{\mathbb{T}^3} f\cdot\Delta u\,dx, 
   \label{Eq2}
\end{align}
where we used $\nabla\cdot u=0$. 
By Lemma~\ref{L1}, the first term on the right side of \eqref{Eq2} is bounded by 
\begin{align*}
     \int_{\mathbb{T}^3}|w| |u| |\Delta u|\,dx
     &\leq
     C\Vert w\Vert_{L^2}^{1/2} \Vert\nabla w\Vert_{L^2}^{1/2} \Vert\nabla u\Vert_{L^2} \Vert \Delta u\Vert_{L^2}
     \leq
     \Vert w\Vert_{L^2} \Vert\nabla w\Vert_{L^2} \Vert\nabla u\Vert_{L^2}^2
     +
     \frac{1}{4} \Vert\Delta u\Vert_{L^2}^2, 
\end{align*} 
where we used Young's inequality. 
By H\"older's inequality, we bound the second term by 
\begin{align*}
     \int_{\mathbb{T}^3}|f| |\Delta u|\,dx
     &\leq
     \Vert f\Vert_{L^2} \Vert \Delta u\Vert_{L^2}
     \leq
     C\Vert f\Vert_{L^2}^2
     +
     \frac{1}{4} \Vert \Delta u\Vert_{L^2}^2.
\end{align*}
Combining all the above estimates and denoting 
$$\bar{X}(t) = \Vert\nabla u(t)\Vert_{L^2}^2 + \alpha^2\Vert \Delta u(t)\Vert_{L^2}^2 \text{\quad and \quad} \bar{Y}(t) = \Vert \Delta u(t)\Vert_{L^2}^2,$$
for $0\leq t\leq T$, 
we arrive at 
\begin{align*}
     \frac{d}{d t}\bar{X}(t) + \bar{Y}(t) 
     \leq
     K_1\bar{X}(t) + C\Vert f\Vert_{L^2}^2,
\end{align*}
where the constant $K_1$ depends on the $H^1$ norms of $u$ and $w$. 
Thus, Gr\"onwall's inequality implies that $u\in L^\infty(0,T;D(A))$.

Next, we show the $H^2$ boundedness of $w$. 
We multiply the $w$ equation in system \eqref{Sys1} by $\Delta^2 w$, integrate by parts over $\mathbb{T}^3$, and obtain 
\begin{align}
     \frac{1}{2}\frac{d}{d t}\left(\Vert\Delta w\Vert_{L^2}^2\right)
     +
     \Vert\nabla\Delta w\Vert_{L^2}^2
     &=
     -
     \int_{\mathbb{T}^3} (u\cdot\nabla)w\cdot \Delta^2 w\,dx
     +
     \int_{\mathbb{T}^3} (w\cdot\nabla)u\cdot \Delta^2 w\,dx \nonumber
     \\&\quad
     +
     \int_{\mathbb{T}^3} (\nabla\times f)\cdot \Delta^2 w\,dx. 
   \label{Eq4}
\end{align}
We then estimate the three terms on the right side of \eqref{Eq4}. 
After integration by parts, 
we use Lemma~\ref{L1} and H\"older's inequality in order to bound the first integral by 
\begin{align*}
     &
     \int_{\mathbb{T}^3} |\nabla u| |\nabla w| |\nabla\Delta w|\,dx
     +
     \int_{\mathbb{T}^3} |u| |\nabla\nabla w| |\nabla\Delta w|\,dx
     \\&\quad
     \leq
     C\Vert\Delta u\Vert_{L^2} \Vert\nabla w\Vert_{L^2}^{1/2}\Vert\Delta w\Vert_{L^2}^{1/2} \Vert\nabla\Delta w\Vert_{L^2}
     +
     C\Vert\nabla u\Vert_{L^2}^{1/2} \Vert\Delta u\Vert_{L^2}^{1/2} \Vert\nabla\nabla w\Vert_{L^2} \Vert\nabla\Delta w\Vert_{L^2}
     \\&\quad
     \leq
     \frac{C}{\sqrt{\lambda_1}}\Vert\Delta u\Vert_{L^2}^2 \Vert\Delta w\Vert_{L^2}^2
     +
     C\Vert\nabla u\Vert_{L^2} \Vert\Delta u\Vert_{L^2} \Vert\Delta w\Vert_{L^2}^2
     +
     \frac{1}{8} \Vert\nabla\Delta w\Vert_{L^2}^2, 
\end{align*} 
where we used \eqref{Agmon} and Young's inequality. 
Similarly, by Lemma~\ref{L1}, the second term is bounded by  
\begin{align*}
     &
     \int_{\mathbb{T}^3} |\nabla w| |\nabla u| |\nabla\Delta w|\,dx
     +
     \int_{\mathbb{T}^3} |w| |\nabla\nabla u| |\nabla\Delta w|\,dx
     \\&\quad
     \leq
     C\Vert\Delta u\Vert_{L^2} \Vert\nabla w\Vert_{L^2}^{1/2} \Vert\Delta w\Vert_{L^2}^{1/2} \Vert\nabla\Delta w\Vert_{L^2}
     +
     C\Vert\nabla\nabla u\Vert_{L^2} \Vert\nabla w\Vert_{L^2}^{1/2} \Vert\Delta w\Vert_{L^2}^{1/2} \Vert\nabla\Delta w\Vert_{L^2}
     \\&\quad
     \leq
     \frac{C}{\sqrt{\lambda_1}}\Vert\Delta u\Vert_{L^2}^2 \Vert\Delta w\Vert_{L^2}^2
     +
     \frac{1}{8} \Vert\nabla\Delta w\Vert_{L^2}^2, 
\end{align*} 
where we also used \eqref{Poincare} and \eqref{Agmon}. 
The last integral is bounded by 
\begin{align*}
     \int_{\mathbb{T}^3}|\nabla(\nabla\times f)| |\nabla\Delta w|\,dx
     &\leq
     \Vert\nabla(\nabla\times f)\Vert_{L^2} \Vert\nabla\Delta w\Vert_{L^2}
     \leq
     C\Vert f\Vert_{H_{\text{curl}}^1}^2
     +
     \frac{1}{8}\Vert\nabla\Delta w\Vert_{L^2}^2,  
\end{align*}
where we integrated by parts and used H\"older's inequality. 
Summing up all the above estimates and denoting 
$$\tilde{X}(t) = \Vert\Delta w(t)\Vert_{L^2}^2 \text{\quad and \quad} \tilde{Y}(t) = \Vert\nabla\Delta w(t)\Vert_{L^2}^2,$$
for $0\leq t\leq T$, we obtain 
\begin{align*}
     \frac{d}{d t}\tilde{X}(t) + \tilde{Y}(t)
     \leq
     K_2\tilde{X}(t) + C\Vert f\Vert_{H_{\text{curl}}^1}^2,
\end{align*}
where the first constant $K_2$ depends on $\Vert u\Vert_{H^2}$ and $\lambda_1$. 
Thus, Gr\"onwall's inequality and the $H^2$ bound of $u$ imply that $w$ is also bounded in $H^2$. 

\subsubsection{\em $H^3$ bounds}
\label{subsubsec4-3-2}
We start by testing \eqref{Sys1_u} with $\Delta^{2}u$, integrating by parts, and obtain 
\begin{align}
     &
     \frac{1}{2}\frac{d}{d t}\left(\Vert \Delta u\Vert_{L^2}^2 + \alpha^2\Vert\nabla \Delta u\Vert_{L^2}^2\right)
     +
     \Vert\nabla \Delta u\Vert_{L^2}^2 
     =
     \int_{\mathbb{T}^3} w\times u\cdot \Delta^{2}u\,dx
     +
     \int_{\mathbb{T}^3} f\cdot \Delta^{2}u\,dx. 
   \label{Eq5}
\end{align}
After integration by parts, the first term on the right side of \eqref{Eq5} is bounded by 
\begin{align*}
     \int_{\mathbb{T}^3}|\nabla w| |u| |\nabla \Delta u|\,dx
     +
     \int_{\mathbb{T}^3}|w| |\nabla u| |\nabla \Delta u|\,dx
     &\leq
     C\Vert u\Vert_{H^1}^{1/2} \Vert u\Vert_{H^2}^{1/2} \Vert\nabla w\Vert_{L^2} \Vert\nabla \Delta u\Vert_{L^2}
     \\&\quad
     +
     C\Vert w\Vert_{H^1}^{1/2} \Vert u\Vert_{H^2}^{1/2} \Vert\nabla u\Vert_{L^2} \Vert\nabla \Delta u\Vert_{L^2}
     \\&
     \leq
     C\Vert u\Vert_{H^1} \Vert u\Vert_{H^2} \Vert\nabla w\Vert_{L^2}^2
     +
     C\Vert w\Vert_{H^1} \Vert w\Vert_{H^2} \Vert\nabla u\Vert_{L^2}^2
     \\&\quad
     +
     \frac{1}{4}\Vert\nabla \Delta u\Vert_{L^2}^2
\end{align*} 
where we applied \eqref{Agmon}. 
By H\"older's inequality and integration by parts, we bound the second term by 
\begin{align*}
     \int_{\mathbb{T}^3}|\nabla f| |\nabla \Delta u|\,dx
     &\leq
     \Vert\nabla f\Vert_{L^2} \Vert\nabla \Delta u\Vert_{L^2}
     \leq
     C\Vert f\Vert_{H^1}^2
     +
     \frac{1}{4} \Vert\nabla \Delta u\Vert_{L^2}^2.
\end{align*} 
Combining the above estimates, we have 
\begin{equation*}
     \frac{d}{dt}\left(\Vert \Delta u\Vert_{L^2}^2 + \alpha^2\Vert\nabla \Delta u\Vert_{L^2}^2\right)
     \leq
     K_3 + C\Vert f\Vert_{H^1}^2, 
\end{equation*}
where the constant $K_3$ depends on the $H^2$ norms of $u$ and $w$. 
Therefore, Gr\"onwall's inequality implies that $u\in L^{\infty}(0, T; H^3\cap V)$. 

Next, we multiply \eqref{Sys1_w} by $\partial^{\beta}w$ after applying the operator $\partial^{\beta}$, 
where $\beta$ is a multi-index with $|\beta|=3$, and obtain 
\begin{align}
     \frac{1}{2}\frac{d}{d t}\left(\Vert\partial^{\beta}w\Vert_{L^2}^2\right)
     +
     \Vert\nabla\partial^{\beta}w\Vert_{L^2}^2
     &=
     -
     \sum_{0<\gamma\leq\beta}{\binom{\beta}{\gamma}}\int_{\mathbb{T}^3} (\partial^{\gamma}u\cdot\nabla)\partial^{\beta-\gamma}w\cdot\partial^{\beta}w\,dx \nonumber
     \\&\quad
     +
     \sum_{0\leq\gamma\leq\beta}{\binom{\beta}{\gamma}}\int_{\mathbb{T}^3} (\partial^{\gamma}w\cdot\nabla)\partial^{\beta-\gamma}u\cdot\partial^{\beta}w\,dx \nonumber
     \\&\quad
     +
     \int_{\mathbb{T}^3} \partial^{\beta}(\nabla\times f)\cdot \partial^{\beta}w\,dx,
   \label{Eq6}
\end{align}
where we used $\nabla\cdot u = 0$ and $\gamma$ is also a multi-index and $\gamma \leq \beta$ indicates that $|\gamma| \leq |\beta|$ and $\gamma_{i}\leq\beta_{i}$ for $i=1,2,3$. 
Then, we estimate the first term on the right side of \eqref{Eq6} in the following three cases. 
For $\gamma \leq \beta$ and $|\gamma|=1$, it is bounded by 
\begin{align*}
     \sum_{|\gamma|=1}{\binom{\beta}{\gamma}}\int_{\mathbb{T}^3} |\partial^{\gamma}u| |\nabla\partial^{\beta-\gamma}w| |\partial^{\beta}w|\,dx
     &
     \leq
     C\Vert\partial^{\gamma}u\Vert_{L^{\infty}} \Vert\nabla\partial^{\beta-\gamma}w\Vert_{L^2} \Vert\partial^{\beta}w\Vert_{L^2}
     \\&
     \leq
     C\Vert u\Vert_{H^3} \Vert w\Vert_{H^{|\beta|}}^2,
\end{align*}
where we used \eqref{Agmon}. 
For $\gamma \leq \beta$ and $|\gamma|=2$, it is bounded by 
\begin{align*}
     \sum_{|\gamma|=2}{\binom{\beta}{\gamma}}\int_{\mathbb{T}^3} |\partial^{\gamma}u| |\nabla\partial^{\beta-\gamma}w| |\partial^{\beta}w|\,dx
     &
     \leq
     C\Vert\partial^{\gamma}u\Vert_{L^{3}} \Vert\nabla\partial^{\beta-\gamma}w\Vert_{L^6} \Vert\partial^{\beta}w\Vert_{L^2}
     \\&
     \leq
     C\Vert u\Vert_{H^2}^{1/2} \Vert u\Vert_{H^3}^{1/2} \Vert w\Vert_{H^{|\beta|}}^2
\end{align*}
For $\gamma \leq \beta$ and $|\gamma|=3$, 
it is bounded from above by 
\begin{align*}
     \sum_{|\gamma|=3}{\binom{\beta}{\gamma}}\int_{\mathbb{T}^3} |\partial^{\gamma}u| |\nabla\partial^{\beta-\gamma}w| |\partial^{\beta}w|\,dx
     &
     \leq
     C\Vert\partial^{\gamma}u\Vert_{L^{2}} \Vert\nabla\partial^{\beta-\gamma}w\Vert_{L^3} \Vert\partial^{\beta}w\Vert_{L^6}
     \\&
     \leq
     C\Vert u\Vert_{H^3} \Vert w\Vert_{H^1}^{1/2} \Vert w\Vert_{H^2}^{1/2} \Vert\nabla\partial^{\beta}w\Vert_{L^2}
     \\&
     \leq
     C\Vert u\Vert_{H^3}^2 \Vert w\Vert_{H^1} \Vert w\Vert_{H^2}
     +
     \frac{1}{8}\Vert\nabla\partial^{\beta}w\Vert_{L^2}^2.
\end{align*}

The estimates for the second term on the right side of \eqref{Eq6} follow similarly in the following four cases. 
For $\gamma \leq \beta$ and $|\gamma|=0$, we integrate by parts and bound it by 
\begin{align*}
     \sum_{|\gamma|=0}{\binom{\beta}{\gamma}}\int_{\mathbb{T}^3} |w| |\partial^{\beta-\gamma}u| |\nabla\partial^{\beta}w|\,dx
     &
     \leq
     C\Vert w\Vert_{L^{\infty}} \Vert\partial^{\beta-\gamma}u\Vert_{L^2} \Vert\nabla\partial^{\beta}w\Vert_{L^2}
     \\&
     \leq
     C\Vert w\Vert_{H^2}^2 \Vert u\Vert_{H^3}^2
     +
     \frac{1}{8}\Vert\nabla\partial^{\beta}w\Vert_{L^2}^2,
\end{align*}
where we used $\nabla\cdot w=0$. 
For $\gamma \leq \beta$ and $|\gamma|=1$, it is bounded by 
\begin{align*}
     \sum_{|\gamma|=1}{\binom{\beta}{\gamma}}\int_{\mathbb{T}^3} |\partial^{\gamma}w| |\nabla\partial^{\beta-\gamma}u| |\partial^{\beta}w|\,dx
     &
     \leq
     C\Vert\partial^{\gamma}w\Vert_{L^{\infty}} \Vert\nabla\partial^{\beta-\gamma}u\Vert_{L^2} \Vert\partial^{\beta}w\Vert_{L^2}
     \\&
     \leq
     C\Vert u\Vert_{H^3} \Vert w\Vert_{H^{|\beta|}}^2, 
\end{align*}
where we applied \eqref{Agmon}. 
For $\gamma \leq \beta$ and $|\gamma|=2$, 
it is estimated in the same way for the first integral with $|\gamma|=2$, i.e., we have 
\begin{align*}
     \sum_{|\gamma|=2}{\binom{\beta}{\gamma}}\int_{\mathbb{T}^3} |\partial^{\gamma}w| |\nabla\partial^{\beta-\gamma}u| |\partial^{\beta}w|\,dx
     &
     \leq
     C\Vert u\Vert_{H^2}^{1/2} \Vert u\Vert_{H^3}^{1/2} \Vert w\Vert_{H^{|\beta|}}^2. 
\end{align*}
For $\gamma \leq \beta$ and $|\gamma|=3$, it is bounded by 
\begin{align*}
     \sum_{|\gamma|=3}{\binom{\beta}{\gamma}}\int_{\mathbb{T}^3} |\partial^{\gamma}w| |\nabla\partial^{\beta-\gamma}u| |\partial^{\beta}w|\,dx
     &
     \leq
     C\Vert\nabla\partial^{\beta-\gamma} u\Vert_{L^{\infty}} \Vert\partial^{\gamma}w\Vert_{L^2} \Vert\partial^{\beta}w\Vert_{L^2}
     \leq
     C\Vert u\Vert_{H^3} \Vert w\Vert_{H^{|\beta|}}^2. 
\end{align*}
As for the last term on the right side of \eqref{Eq6}, we integrate by parts and estimate as 
\begin{align*}
     \int_{\mathbb{T}^3} \partial^{\beta}(\nabla\times f)\cdot \partial^{\beta}w\,dx
     &
     \leq
     \sum_{\gamma=2}{\binom{\beta}{\gamma}}\int_{\mathbb{T}^3}|\partial^{\gamma}(\nabla\times f)| |\partial\partial^{\beta}w|\,dx
     \leq
     C\Vert\nabla\times f\Vert_{H^2} \Vert\nabla\partial^{\beta}w\Vert_{L^2}
     \\&
     \leq
     C\Vert f\Vert_{H_{\text{curl}}^2}^2
     +
     \frac{1}{8}\Vert\nabla\partial^{\beta}w\Vert_{L^2}^2.
\end{align*}
Summing up all the above estimates, we have 
\begin{equation*}
     \frac{d}{dt}\Vert\partial^{\beta}w\Vert_{L^2}^2 + \Vert\nabla\partial^{\beta}w\Vert_{L^2}^2
     \leq
     K_4\Vert\partial^{\beta}u\Vert_{L^2}^2 + K_5, 
\end{equation*}
where the constants $K_4$ and $K_5$ depend on the $H^3$ norm of $u$, $H^2$ norm of $w$, 
while $K_5$ also depends on the $H_{\text{curl}}^2$ norm of $f$. 
Therefore, Gr\"onwall's inequality implies that $w\in L^{\infty}(0, T; H^3\cap V)\cap L^{2}(0, T; H^4\cap V)$. 
Therefore, by repeating similar arguments as above inductively, we get $H^s$ uniform bound on $u$ and $w$ for all integers $s\geq 4$.  
Proof of Theorem~\ref{T2} is thus complete. 

\section{Proof of Convergence Results}
\label{sec5}
In this section, we prove our convergence results Theorem~\ref{T3} and Theorem~\ref{T4}. 
Notice that from the proof of Theorem~\ref{T1}, we have $\nabla\cdot w(t) = 0$ for a.e. $t>0$ as long as $\nabla\cdot w_0= 0$. 

\subsection{Convergence of $\nabla\times u$ to $w$ in $L^2$}
\label{subsec5-1}
{\smallskip\noindent {\em Proof of Theorem~\ref{T3}.}} 

We start by applying the {\it{curl}} operator ``$\nabla\times$" to  \eqref{Sys1_u} and after denoting by $\omega$ the vorticity $\nabla\times u$, 
we obtain 
\begin{align}
     &
     (I - \alpha^2\Delta)\omega_{t} 
     -
     \nu\Delta \omega
     + 
     (u\cdot\nabla)w
     -
     (\nabla\cdot w)u
     -
     (w\cdot\nabla)u 
     = 
     \nabla\times f,
   \label{Curl}
\end{align}
where we used $\nabla\cdot u = 0$ and the identity 
$$\nabla\times({\bf{F}}\times{\bf{G}}) = \left( (\nabla\cdot{\bf{G}}) + {\bf{G}}\cdot\nabla \right){\bf{F}} - \left( (\nabla\cdot{\bf{F}}) + {\bf{F}}\cdot\nabla \right){\bf{G}}$$
for arbitrary smooth vector fields ${\bf{F}}$ and ${\bf{G}}$ in $\mathbb{R}^3$. 
Denoting by $\xi$ the difference $\omega - w$ and subtracting the $w$ equation of system \eqref{Sys1} from \eqref{Curl} lead to 
\begin{align}
     &
     \xi_{t}
     -
     \alpha^2\Delta\omega_{t}
     -
     \Delta\xi
     -
     (\nabla\cdot w)u
     =
     0. 
   \label{Curl_diff}
\end{align}
Then, we use Theorem~\ref{T1} and rewrite \eqref{Curl_diff} as 
 \begin{align*}
     &
     \xi_{t}
     -
     \alpha^2\Delta\xi_{t}
     -
     \Delta\xi
     =
     \alpha^2\Delta w_{t}
     +
     (\nabla\cdot w)u, 
\end{align*}
to which we multiply $\xi$ and integrate by parts over $\mathbb{T}^3$, and obtain 
 \begin{align}
     &
     \frac{1}{2}\frac{d}{d t}\left(\Vert\xi\Vert_{L^2}^2 + \alpha^2\Vert\nabla\xi\Vert_{L^2}^2\right)
     +
     \Vert\nabla\xi\Vert_{L^2}^2
     =
     \alpha^2\int_{\mathbb{T}^3} \Delta w_{t}\cdot\xi\,dx
     +
     \int_{\mathbb{T}^3} (\nabla\cdot w)u\cdot\xi\,dx.  
   \label{Eq7}
\end{align}
Note that the second term on the right side of \eqref{Eq7} vanishes due to Theorem \ref{T1} since $\nabla\cdot w(0) = 0$. 
In order to estimate the first term on the right side of \eqref{Eq7}, 
we integrate by parts and use the equation of $w$ and obtain 
\begin{align}
     \alpha^2\int_{\mathbb{T}^3} \Delta w_{t}\cdot\xi\,dx
     =
     \alpha^2\int_{\mathbb{T}^3} w_{t}\cdot\Delta\xi\,dx
     &=
     \alpha^2\int_{\mathbb{T}^3}\Delta w\cdot\Delta\xi\,dx
     -
     \alpha^2\int_{\mathbb{T}^3}(u\cdot\nabla)w\cdot\Delta\xi\,dx \nonumber
     \\&\quad\quad
     +
     \alpha^2\int_{\mathbb{T}^3}(w\cdot\nabla)u\cdot\Delta\xi\,dx
     +
     \alpha^2\int_{\mathbb{T}^3}\nabla\times f\cdot\Delta\xi\,dx.
   \label{Eq8}
\end{align}
Then, we estimates the four integrals on the right side of \eqref{Eq8}. 
After integration by parts, we bound the first integral by 
\begin{align*}
     \alpha^2\int_{\mathbb{T}^3} |\nabla\Delta w|\,|\nabla\xi|\,dx
     &\leq
     C\alpha^2\Vert\nabla\Delta w\Vert_{L^2} \Vert\nabla\xi\Vert_{L^2}
     \leq
     C\alpha^2\Vert w\Vert_{H^3}^2
     +
     \alpha^2\Vert\nabla\xi\Vert_{L^2}^2.
\end{align*}
Using Lemma~\ref{L1}, the second integral is bounded by 
\begin{align*}
     &
     \alpha^2\int_{\mathbb{T}^3} |\nabla u|\,|\nabla w|\,|\nabla\xi|\,dx
     +
     \alpha^2\int_{\mathbb{T}^3} |u|\,|\Delta w|\,|\nabla\xi|\,dx
     \\&\quad
     \leq
     C\alpha^2\Vert\Delta u\Vert_{L^2} \Vert w\Vert_{H^2} \Vert\nabla\xi\Vert_{L^2}
     +
     C\alpha^2\Vert\Delta u\Vert_{L^2} \Vert\Delta w\Vert_{L^2} \Vert\nabla\xi\Vert_{L^2}
     \\&\quad
     \leq
     C\alpha^2\Vert u\Vert_{H^2}^2 \Vert w\Vert_{H^2}^2
     +
     \alpha^2\Vert\nabla\xi\Vert_{L^2}^2.
\end{align*}
Estimates for the third integral is similar and we have 
\begin{align*}
     \alpha^2\int_{\mathbb{T}^3}(w\cdot\nabla)u\,\Delta\xi\,dx
     &
     \leq
     C\alpha^2\Vert\Delta u\Vert_{L^2} \Vert w\Vert_{H^2} \Vert\nabla\xi\Vert_{L^2}
     +
     C\alpha^2\Vert\Delta u\Vert_{L^2} \Vert w\Vert_{H^2} \Vert\nabla\xi\Vert_{L^2}
     \\&
     \leq
     C\alpha^2\Vert u\Vert_{H^2}^2 \Vert w\Vert_{H^2}^2
     +
     \alpha^2\Vert\nabla\xi\Vert_{L^2}^2.
\end{align*}
As for the last integral in \eqref{Eq8}, 
we integrate by parts and use H\"older's inequality, and bound it by 
\begin{align*}
     \alpha^2\int_{\mathbb{T}^3}|\Delta f|\,|\nabla\xi|\,dx
     &
     \leq
     C\alpha^2\Vert\Delta f\Vert_{L^2} \Vert\nabla\xi\Vert_{L^2}
     \leq
     C\alpha^2\Vert f\Vert_{H^2}^2
     +
     \alpha^2\Vert\nabla\xi\Vert_{L^2}^2. 
\end{align*}
Summing up all the above estimates, we obtain 
\begin{align}
     &
     \frac{d}{d t}\left(\Vert\xi\Vert_{L^2}^2 + \alpha^2\Vert\nabla\xi\Vert_{L^2}^2\right)
     \leq
     C\left(\Vert\xi\Vert_{L^2}^2 + \alpha^2\Vert\nabla\xi\Vert_{L^2}^2\right)
     +
     \tilde{K}\alpha^2
  \label{Eq9}
\end{align}
where the constant $\tilde{K}$ depends on the $H^3$ norms of $u$ and $w$, as well as $H^2$ norm of $f$. 
By Lemma~\ref{L3} we have 
\begin{align*}
     &
     \Vert\xi(t)\Vert_{L^2}^2 + \alpha^2\Vert\nabla\xi(t)\Vert_{L^2}^2
     + 
     \int_0^t\Vert\nabla\xi(s)\Vert_{L^2}^2\,ds
     \leq
     \tilde{K}\alpha^2(e^{Ct} - 1), 
\end{align*}
where we used $\xi_0 = w_0 - \omega_0 = 0$. 
Therefore, 
$\Vert\xi(t)\Vert_{L^\infty(0,T;H)}^2+
\Vert\xi(t)\Vert_{L^2(0,T;V)}^2 \leq C\alpha^2 e^{CT} \to 0$ as $\alpha\to 0$.
The proof of Theorem~\ref{T3} is thus complete. 

\subsection{Convergence of $\omega$ to $\tilde{\omega}$ and $u$ to $\tilde{u}$ in $L^2$}
\label{subsec5-2}
{\smallskip\noindent {\em Proof of Theorem~\ref{T4}.}} 
Assume the hypotheses and the notation of Theorem \ref{T4}.  Since $u_0\in H^4\cap V$, by Theorem \ref{thm_NSE_local_reg}, there exists a time $T>0$ and a unique strong solution $(\tilde u,\tilde p)$ to \eqref{NSE} satisfying $\tilde u \in C([0,T];H^4\cap V)\cap  L^2(0,T;H^{5}\cap V)$.  

%\subsubsection{Vorticity convergence}
%\label{subsubsec5-2-1}
In view of Theorem~\ref{T3}, it suffices to show that $\Vert w - \tilde{\omega}\Vert_{L^2} + \Vert u - \tilde{u}\Vert_{L^2} \sim\mathcal{O}(\alpha)$. 
We start by applying the curl operator ``$\nabla\times$" to \eqref{NSE_u} and obtain  
\begin{align}
     &
     \frac{\partial \tilde{\omega}}{\partial t} 
     -
     \Delta \tilde{\omega}
     + 
     (\tilde{u}\cdot\nabla)\tilde{\omega}
     -
     (\tilde{\omega}\cdot\nabla)\tilde{u} 
     = 
     \nabla\times f. 
   \label{vorticity}
\end{align}
Then, by taking the difference of \eqref{Sys1_w} and \eqref{vorticity} and denoting by $\theta:=w - \tilde{\omega}$, we have 
\begin{align*}
     &
     \frac{\partial\theta}{\partial t}
     -
     \Delta\theta
     +
     (u\cdot\nabla)w
     -
     (\tilde{u}\cdot\nabla)\tilde{\omega}
     -
     (w\cdot\nabla)u
     +
     (\tilde{\omega}\cdot\nabla)\tilde{u}
     =
     0.
\end{align*}
Denoting $\zeta:=u-\tilde{u}$, we rewrite the above as
\begin{align}
     &
     \frac{\partial\theta}{\partial t}
     -
     \Delta\theta
     +
     (u\cdot\nabla)\theta
     +
     (\zeta\cdot\nabla)\tilde{\omega}
     -
     (\theta\cdot\nabla)u
     -
     (\tilde{\omega}\cdot\nabla)\zeta
     =
     0, 
     \label{vor_diff}
\end{align}
with $\theta(\cdot, 0) = \theta_0 = 0$. 

%\subsubsection{Velocity convergence}
%\label{subsubsec5-2-2}
Next, by subtracting \eqref{NSE_vor} from \eqref{Sys1_u}, we obtain the following system for $\zeta:=u-\tilde{u}$.
\begin{subequations}
\begin{empheq}[left=\empheqlbrace]{align}  
     &
     (I - \alpha^2\Delta)\frac{\partial \zeta}{\partial t}
     -
     \alpha^2\Delta\frac{\partial\tilde{u}}{\partial t}
     -
     \Delta\zeta
     +
     w\times u - \tilde{\omega}\times\tilde{u}
     +
     \nabla\Pi
     =
     0, \label{u_diff}
     \\&
     \nabla\cdot\zeta = 0, \label{u_diff_div_free}
     \\&
     \zeta(\cdot, 0) = \zeta_0 = 0, \label{u_diff_ini}
\end{empheq}
\end{subequations}
where $\Pi = p - \tilde{p} - \frac{1}{2}|\tilde{u}|^2$. 
We multiply \eqref{vor_diff} by $\theta$ and \eqref{u_diff} by $\zeta$, respectively, 
integrate by parts over $\mathbb{T}^3$, and add, to obtain
\begin{align}
     &
     \frac{1}{2}\frac{d}{dt}\left(\Vert\theta\Vert_{L^2}^2 + \Vert\zeta\Vert_{L^2}^2 + \alpha^2\Vert\nabla\zeta\Vert_{L^2}^2\right)
     +
     \Vert\nabla\zeta\Vert_{L^2}^2
     +
     \Vert\nabla\theta\Vert_{L^2}^2 \nonumber
     \\&
     =
     -\int_{\mathbb{T}^3}(\zeta\cdot\nabla)\tilde{\omega}\cdot\theta\,dx
     +
     \int_{\mathbb{T}^3}(\theta\cdot\nabla)u\cdot\theta\,dx
     +
     \int_{\mathbb{T}^3}(\tilde{\omega}\cdot\nabla)\zeta\cdot\theta\,dx \nonumber
     \\&\quad
     +
     \alpha^2\int_{\mathbb{T}^3}\Delta\frac{\partial\tilde{u}}{\partial t}\cdot\zeta\,dx 
     +
     \int_{\mathbb{T}^3}(\theta\times \tilde{u})\cdot\zeta\,dx, 
     \label{u_diff_energy}
\end{align}
where we used $w\times u - \tilde{\omega}\times\tilde{u} = \theta\times \tilde{u} + w\times\zeta$, 
$(w\times\zeta)\cdot\zeta = 0$, \eqref{embd2}, and \eqref{u_diff_div_free}. 
Next, we estimate the five integrals on the right side of \eqref{u_diff_energy}. 
Using \eqref{B:326}, the first integral is bounded by 
\begin{align*}
%      \int_{\mathbb{T}^3}|\tilde{\omega}| |\zeta| |\nabla\theta|\,dx
%      &\leq
     \Vert\nabla\tilde{\omega}\Vert_{L^2} \Vert\zeta\Vert_{L^2}^{1/2} \Vert\nabla\zeta\Vert_{L^2}^{1/2} \Vert\nabla\theta\Vert_{L^2}
     \leq
     C\Vert\zeta\Vert_{L^2}^{2}
     +
     \frac{1}{4}\Vert\nabla\theta\Vert_{L^2}^2
     +
     \frac{1}{6}\Vert\nabla\zeta\Vert_{L^2}^2,
\end{align*}
where $\Vert\nabla\tilde{\omega}\Vert_{L^\infty(0,T;H)}\leq C$.  
The second integral can also be estimated using \eqref{B:623s}:
\begin{align*}
     \int_{\mathbb{T}^3}(\theta\cdot\nabla)u\cdot\theta\,dx
     &\leq
     \Vert\nabla\theta\Vert_{L^2}^{2} \Vert u\Vert_{L^2}^{1/2} \Vert\nabla u\Vert_{L^2}^{1/2} 
     \leq
     C\Vert \nabla u\Vert_{L^2}\Vert\theta\Vert_{L^2}^{2}.
\end{align*}
We note that, as in Remark \ref{remark_L2V_bdd}, that $\Vert \nabla u\Vert_{L^2}$ might not be bounded independently of $\alpha$, but $\int_0^T\Vert \nabla u(t)\Vert_{L^2}\,dt\leq T^{1/2}\|u\|_{L^2(0,T;V)}$ is bounded independently of $\alpha\in(0,1]$.  
The third integral can be estimated using \eqref{B:623}:
\begin{align*}
     \int_{\mathbb{T}^3}(\tilde{\omega}\cdot\nabla)\zeta\cdot\theta\,dx
     &\leq
     \Vert\nabla\tilde{\omega}\Vert_{L^2} \Vert\theta\Vert_{L^2}^{1/2} \Vert\nabla\theta\Vert_{L^2}^{1/2} \Vert\nabla\zeta\Vert_{L^2}
     \leq
     C\Vert\theta\Vert_{L^2}^{2}
     +
     \frac{1}{4}\Vert\nabla\theta\Vert_{L^2}^2
     +
     \frac{1}{6}\Vert\nabla\zeta\Vert_{L^2}^2.
\end{align*}
By substituting ${\partial \tilde{u}}/{\partial t}$ from \eqref{NSE_u} and integration by parts, 
we bound the fourth term on the right side of \eqref{u_diff_energy} by
\begin{align*}
     &
     \alpha^2\int_{\mathbb{T}^3}|\nabla\Delta\tilde{u}||\nabla\zeta|\,dx
     +
     \alpha^2\int_{\mathbb{T}^3}|\nabla(\tilde{\omega}\times\tilde{u})||\nabla\zeta|\,dx
     +
     \alpha^2\int_{\mathbb{T}^3}|\nabla f||\nabla\zeta|\,dx
     \\&
     \leq
     C\alpha^2\Vert\tilde{u}\Vert_{H^3}^2
     +
     C\alpha^2\Vert\tilde{u}\Vert_{H^1}^2\Vert\tilde{u}\Vert_{H^2}\Vert\tilde{u}\Vert_{H^3}
     +
     C\alpha^2\Vert\tilde{u}\Vert_{H^1}\Vert\tilde{u}\Vert_{H^2}^3
     +
     C\alpha^2\Vert f\Vert_{H^1}^2
     +
     \frac{1}{6}\Vert\nabla\zeta\Vert_{L^2}^2,
\end{align*} 
where we used the hypothesis of the theorem that $\alpha\in(0,1]$.
As for the last integral, we use \eqref{Agmon} to bound it by 
\begin{align*}
     \Vert \tilde{u}\Vert_{H^2} \Vert\theta\Vert_{L^2} \Vert\zeta\Vert_{L^2}
     \leq
     C\Vert\theta\Vert_{L^2}^2
     +
     C\Vert\zeta\Vert_{L^2}^2,
\end{align*}
where $\Vert \tilde{u}\Vert_{H^2}\leq C$.  
Combining all the above estimates and using \eqref{u_diff_ini}, we obtain 
\begin{align*}
&\quad
\frac{d}{dt}\left(\Vert\theta(t)\Vert_{L^2}^2 + \Vert\zeta(t)\Vert_{L^2}^2 + \alpha^2\Vert\nabla\zeta(t)\Vert_{L^2}^2\right) 
+ \Vert\nabla\zeta\Vert_{L^2}^2 
+ \Vert\nabla\theta\Vert_{L^2}^2 
\\&\leq 
C\alpha^2 
+ 
\tilde{C}\Vert \nabla u\Vert_{L^2}\left(\Vert\theta(t)\Vert_{L^2}^2 
+ \Vert\zeta(t)\Vert_{L^2}^2 
+ \alpha^2\Vert\nabla\zeta(t)\Vert_{L^2}^2\right),
\end{align*}
where the constants $C$ and $\tilde{C}$ are uniform-in-time bounds for $\Vert\tilde{u}\Vert_{H^3}$, as well as $\tilde{K}$. 
Using Gr\"onwall's inequality, the integrability of $\Vert \nabla u(t)\Vert_{L^2}$, and the fact that $\|u\|_{L^2(0,T;V)}$ is bounded independently of $\alpha\in(0,1]$, we obtain that $\Vert\theta(t)\Vert_{L^2}^2 + \Vert\zeta(t)\Vert_{L^2}^2 
+ \int_0^t\Vert\nabla\zeta(s)\Vert_{L^2}^2\,ds
+ \int_0^t\Vert\nabla\theta(s)\Vert_{L^2}^2\,ds
\leq 
C\alpha^2T$. 
Thus, the proof of Theorem~\ref{T4} is now complete.  

\subsection{Proof of Corollary~\ref{blow_up}}
Assume the hypothses.  From Theorem \ref{T4}, $\Vert u - \tilde{u}\Vert_{L^\infty(0,T;H)} \to 0$ and $\Vert u - \tilde{u}\Vert_{L^2(0,T;V)}$ as $\alpha \to 0$.  Thus, passing to the limit as $\alpha\to0$ in the energy equality \eqref{energy_equality_u} yields
\begin{align}\label{blow_up_lim}
  \limsup_{\alpha\rightarrow 0}\alpha^2\|\nabla u(t)\|_{L^2}^2 
 + \|\tilde u(t)\|_{L^2}^2 
 +2\int_0^t\|\nabla \tilde u(s)\|_{L^2}^2\,ds
 =
 \|u_0\|_{L^2}^2+
 2\int_0^t( \tilde u(s),f(s))\,ds
\end{align}
However, if $\tilde u$ is a strong solution to the Navier-Stokes equation on $[0,T]$, the energy identity
\begin{align}\label{energy_equality_NSE}
  \|\tilde u(t)\|_{L^2}^2 
 +2\int_0^t\|\nabla \tilde u(s)\|_{L^2}^2\,ds
 =
 \|u_0\|_{L^2}^2
 +2\int_0^t( \tilde u(s),f(s))\,ds,
\end{align}
holds (see, e.g., \cite{Constantin_Foias_1988, Foias_Manley_Rosa_Temam_2001}).  Thus, \eqref{blow_up_lim} together with \eqref{energy_equality_NSE} contradict \eqref{blowup_criterion}.

\section*{Acknowledgments} 

\noindent Author A.L. is partially supported by NSF grant DMS-1716801. 

\noindent Author L.R. is partially supported by NSF grant DMS-1522191.
\bibliographystyle{abbrv}

\end{document}